\documentclass[reqno,centertags,draft]{amsart}
\usepackage{amsmath,amsthm,amscd,amssymb,latexsym,upref,enumerate}
\newcommand{\bbC}{{\mathbb C}}
\newcommand{\bbD}{{\mathbb D}}

\newcommand{\bbF}{{\mathbb F}}

\newcommand{\bbN}{{\mathbb N}}

\newcommand{\bbR}{{\mathbb R}}

\newcommand{\bbZ}{{\mathbb Z}}

\newcommand{\cE}{{\mathcal E}}

\newcommand{\cM}{{\mathcal M}}
\newcommand{\cN}{{\mathcal N}}
\newcommand{\cO}{{\mathcal O}}

\newcommand{\cR}{{\mathcal R}}

\newcommand{\rank}{\text{\rm{rank}}}

\newcommand{\Arg}{\text{\rm{Arg}}}
\newcommand{\Arc}{\text{\rm{Arc}}}
\newcommand{\dom}{\text{\rm{dom}}}
\newcommand{\tr}{\text{\rm{tr}}}

\newcommand{\ess}{\text{\rm{ess}}}
\newcommand{\ac}{\text{\rm{ac}}}
\renewcommand{\sc}{\text{\rm{sc}}}
\newcommand{\pp}{\text{\rm{pp}}}
\newcommand{\disc}{\text{\rm{disc}}}
\newcommand{\supp}{\text{\rm{supp}}}
\newcommand{\dist}{\text{\rm{dist}}}
\newcommand{\diam}{\text{\rm{diam}}}

\newcommand{\id}{\text{\rm{id}}}
\newcommand{\ca}{\text{\rm{cap}}}

\newcommand{\Om}{\Omega}
\newcommand{\om}{\omega}
\newcommand{\si}{\sigma}
\newcommand{\la}{\lambda}
\newcommand{\La}{\Lambda}
\newcommand{\al}{\alpha}
\newcommand{\be}{\beta}
\newcommand{\Ga}{\Gamma}
\newcommand{\ga}{\gamma}

\newcommand{\de}{\delta}
\newcommand{\te}{\theta}
\newcommand{\ze}{\zeta}

\newcommand{\eps}{\varepsilon}

\newcommand{\Oh}{O}
\newcommand{\oh}{o}

\newcommand{\bi}{\bibitem}
\newcommand{\no}{\notag}
\newcommand{\lb}{\label}
\newcommand{\f}{\frac}

\newcommand{\ol}{\overline}

\newcommand{\wti}{\widetilde}
\newcommand{\bs}{\backslash}

\newcommand{\abs}[1]{\left\lvert#1\right\rvert}

\newcommand{\dott}{\,\cdot\,}
\newcommand{\dD}{{\partial\mathbb{D}}}
\newcommand{\dF}{{\partial\mathbb{F}}}
\newcommand{\dOm}{{\partial\Omega}}
\newcommand{\Cl}{\mathbb{C}_{\ell}}
\newcommand{\Cr}{\mathbb{C}_{r}}

\newcommand{\st}{\;|\;}
\newcommand{\ltz}{{\ell^2(\bbZ)}}
\newcommand{\lt}[1]{{\ell^2(#1)}}
\newcommand{\cm}{{\bf{x}}}

\newcommand{\deven}{\delta_{\rm even}}
\newcommand{\dodd}{\delta_{\rm odd}}

\renewcommand{\Re}{\text{\rm Re}}
\renewcommand{\Im}{\text{\rm Im}}
\renewcommand{\ln}{\text{\rm ln}}

\allowdisplaybreaks
\numberwithin{equation}{section}

\newtheorem{theorem}{Theorem}[section]

\newtheorem{lemma}[theorem]{Lemma}

\newtheorem{hypothesis}[theorem]{Hypothesis}
\newtheorem{definition}[theorem]{Definition}

\theoremstyle{definition}

\begin{document}

\title[Local Spectral Properties of reflectionless operators]
{Local Spectral Properties of Reflectionless Jacobi, CMV, \\ and
Schr\"odinger Operators}
\author[F.\ Gesztesy and M.\ Zinchenko]
{Fritz Gesztesy and Maxim Zinchenko}

\address{Department of Mathematics,
University of Missouri, Columbia, MO 65211, USA}
\email{fritz@math.missouri.edu}
\urladdr{http://www.math.missouri.edu/personnel/faculty/gesztesyf.html}

\address{Department of Mathematics,
California Institute of Technology, Pasadena, CA 91125, USA}
\email{maxim@caltech.edu}
\urladdr{http://math.caltech.edu/\~{}maxim}

%\thanks{}
\date{\today}
\dedicatory{Dedicated with great pleasure to Ludwig Streit on the occasion of his 70th birthday.}
\thanks{To appear in {\it J. Diff. Eq.}}

%%%%%%%%%%%%%%%%%%%%%%%%%%%%%%%%%%%%%%%
%%%%%%%%%%%%%%%%%%%%%%%%%%%%%%%%%%%%%%%
\begin{abstract}
We prove that Jacobi, CMV, and Schr\"odinger operators, which are
reflectionless on a homogeneous set $\cE$ (in the sense of Carleson), under the assumption of a Blaschke-type condition on their discrete spectra accumulating
at $\cE$, have purely absolutely continuous spectrum on $\cE$.
\end{abstract}
%%%%%%%%%%%%%%%%%%%%%%%%%%%%%%%%%%%%%%%
%%%%%%%%%%%%%%%%%%%%%%%%%%%%%%%%%%%%%%%

\maketitle

%%%%%%%%%%%%%%%%%%%%%%%%%%%%%%%%%%%%%%%
%%%%%%%%%%%%%%%%%%%%%%%%%%%%%%%%%%%%%%%
\section{Introduction}\lb{s1}
%%%%%%%%%%%%%%%%%%%%%%%%%%%%%%%%%%%%%%%
%%%%%%%%%%%%%%%%%%%%%%%%%%%%%%%%%%%%%%%

In this paper we consider self-adjoint Jacobi and Schr\"odinger operators $H$
on $\bbZ$ and $\bbR$, respectively, and unitary CMV operators $U$ on $\bbZ$,
which are reflectionless on a homogeneous set $\cE$ contained in the essential spectrum. We prove that under the assumption of a Blaschke-type
condition on their discrete spectra accumulating at $\cE$, the operators $H$,
respectively, $U$, have purely absolutely continuous spectrum on $\cE$.

We note that homogeneous sets were originally discussed by Carleson \cite{Ca83};
we also refer to \cite{JM85}, \cite{PY03}, and \cite{Zi89} in this context. Morever,
by results of Kotani \cite{Ko84}--\cite{Ko87b} (also recorded in detail in
\cite[Theorem 12.5]{PF92}), it is known that CMV, Jacobi, and Schr\"odinger operators, reflectionless on a set $\cE$ of positive Lebesgue measure, have absolutely continuous spectrum on the essential closure of $\cE$, denoted by
$\ol{\cE}^e$ (with uniform multiplicity two on $\cE$). This result has recently been revisited in \cite{GMZ08}. The focal point of this paper is to show that under suitable additional conditions on $\cE$, such
as $\cE$ homogeneous, and a Blaschke-type condition on the discrete spectrum accumulating at $\cE$, the spectrum is actually purely absolutely continuous on $\cE$.

To put this result in some perspective, we briefly single out Schr\"odinger operators and illustrate the notion of being reflectionless: {\it Reflectionless} (self-adjoint)
Schr\"odinger operators $H$ in $L^2(\bbR;dx)$ can be characterized, for instance, by the fact that for all $x\in\bbR$ and for a.e.\ $\lambda\in \sigma_{\rm ess}(H)$, the diagonal Green's function of $H$ has purely imaginary normal boundary values,
\begin{equation}
G(\lambda+i0,x,x) \in i\bbR.    \lb{1.3}
\end{equation}
Here $\sigma_{\rm ess}(H)$ denotes the essential spectrum of $H$ (we assume
$\sigma_{\rm ess}(H)\neq\emptyset$) and
\begin{equation}
G(z,x,x')=(H-zI)^{-1}(x,x'), \quad z\in\bbC\backslash\sigma(H), \lb{1.4}
\end{equation}
denotes the integral kernel of the resolvent of $H$. This global notion of reflectionless Schr\"odinger operators can of course be localized and extends to subsets of
$\sigma_{\rm ess}(H)$ of positive Lebesgue measure. In the actual body of our paper we will use an alternative definition of the notion of reflectionless Schr\"odinger operators conveniently formulated directly in terms of half-line Weyl--Titchmarsh functions; we refer to Definitions \ref{d2.2}, \ref{d3.2}, and \ref{d4.2} for more details. For various discussions of classes of reflectionless differential and difference operators, we refer, for instance, to Craig \cite{Cr89}, De Concini and Johnson \cite{DJ87},
Deift and Simon \cite{DS83},
Gesztesy, Krishna, and Teschl \cite{GKT96}, Gesztesy and Yuditskii \cite{GY06},
Johnson \cite{Jo82}, Kotani \cite{Ko84}, \cite{Ko87a}, Kotani and Krishna \cite{KK88}, Peherstorfer and Yuditskii \cite{PY06}, Remling \cite{Re07}, \cite{Re08}, Sims \cite{Si07a}, and Sodin and Yuditskii \cite{SY95a}--\cite{SY96}. In particular, we draw attention to the recent papers by Remling \cite{Re07}, \cite{Re08}, that illustrate in great depth the ramifications of the existence of absolutely continuous spectra in one-dimensional problems.

The trivial case $H_0=-d^2/dx^2$, and the $N$-soliton potentials $V_N$, $N\in\bbN$,
that is, exponentially decreasing solutions in $C^\infty(\bbR)$ of some (and
hence infinitely many) equations of the stationary Korteweg--de Vries (KdV)
hierarchy, yield well-known examples of reflectionless Schr\"odinger
operators $H_N=-d^2/dx^2+V_N$. Similarly, all periodic Schr\"odinger
operators are reflectionless. Indeed, if $V_a$ is periodic with some period $a>0$,
that is, $V_a(x+a)=V_a(x)$ for a.e.\ $x\in\bbR$, then standard Floquet
theoretic considerations show that the spectrum of $H_a=-d^2/dx^2+V_a$ is a
countable union of compact intervals (which may degenerate into a union of
finitely-many compact intervals and a half-line) and the diagonal Green's function
of $H_a$ is purely imaginary for every point in the open interior of $\sigma(H_a)$. More
generally, certain classes of quasi-periodic and almost periodic
potentials also give rise to reflectionless Schr\"odinger operators with
homogeneous spectra. The prime example of such quasi-periodic potentials is
represented by the class of real-valued bounded algebro-geometric KdV
potentials corresponding to an underlying (compact) hyperelliptic Riemann
surface (see, e.g., \cite[Ch.\ 3]{BBEIM94}, \cite{DMN76},
\cite[Ch.\ 1]{GH03}, \cite{Jo88}, \cite[Chs.\ 8, 10]{Le87}, \cite[Ch.\ 4]{Ma86},
\cite[Ch.\ II]{NMPZ84} and the literature cited therein). These examples yield reflectionless operators in a global sense, that is, they are reflectionless on the whole spectrum. On the other hand, as discussed, recently by Remling
\cite{Re07}, the notion of being reflectionless also makes sense locally on subsets of the spectrum. More general classes of almost periodic Schr\"odinger operators, reflectionless on sets where the Lyapunov exponent vanishes, were studied by  Avron and Simon \cite{AS81}, Carmona and Lacroix \cite[Ch.\ VII]{CL90},
Chulaevskii \cite{Ch84}, Craig \cite{Cr89},
Deift and Simon \cite{DS83},  Egorova \cite{Eg92}, Johnson \cite{Jo82},
Johnson and Moser \cite{JM82}, Kotani \cite{Ko84}--\cite{Ko87b},
Kotani and Krishna \cite{KK88}, Levitan \cite{Le82}--\cite{Le85},
\cite[Chs.\ 9, 11]{Le87},  Levitan and Savin \cite{LS84}, Moser \cite{Mo81},
Pastur and Figotin \cite[Chs.\ V, VII]{PF92}, Pastur and Tkachenko
\cite{PT89}, and Sodin and Yuditskii \cite{SY95a}--\cite{SY96}.

Analogous considerations apply to Jacobi operators (see, e.g., \cite{CL90},
\cite{Te00} and the literature cited therein) and CMV operators (see
\cite{Si04a}--\cite{Si07} and the extensive list of references provided therein and \cite{GZ06a} for the notion of reflectionless CMV operators).

In Section \ref{s2} we consider the case of Jacobi operators; CMV operators are studied in Section \ref{s3} followed by Schr\"odinger operators in Section \ref{s4}. Herglotz and Weyl--Titchmarsh functions in connection with Jacobi and Schr\"odinger Operators are discussed in Appendix \ref{B}; Caratheodory and Weyl--Titchmarsh functions for CMV operators are summarized in Appendix \ref{C}.

%%%%%%%%%%%%%%%%%%%%%%%%%%%%%%%%%%%%%%%
%%%%%%%%%%%%%%%%%%%%%%%%%%%%%%%%%%%%%%%
\section{Reflectionless Jacobi Operators} \lb{s2}
%%%%%%%%%%%%%%%%%%%%%%%%%%%%%%%%%%%%%%%
%%%%%%%%%%%%%%%%%%%%%%%%%%%%%%%%%%%%%%%

In this section we investigate spectral properties of self-adjoint
Jacobi operators reflectionless on compact homogeneous subsets of
the real line.

We start with some general considerations of self-adjoint Jacobi
operators. Let $a=\{a(n)\}_{n\in\bbZ}$ and $b=\{b(n)\}_{n\in\bbZ}$
be two sequences (Jacobi parameters) satisfying
\begin{equation}
a, b \in \ell^\infty(\bbZ), \quad a(n)>0, \; b(n)\in\bbR, \;
n\in\bbZ,
\end{equation}
and denote by $L$ the second-order difference expression defined by
\begin{equation}
L = a S^+ + a^- S^- + b,
\end{equation}
where we use the notation for $f=\{f(n)\}_{n\in\bbZ}\in\ell^{\infty}(\bbZ)$,
\begin{equation}
(S^\pm f)(n)=f(n\pm 1)= f^\pm(n), \; n\in\bbZ, \quad
S^{++}=(S^+)^+, \; S^{--}=(S^-)^-, \, \text{ etc.}   \lb{2.3}
\end{equation}
Moreover, we introduce the associated bounded self-adjoint Jacobi
operator $H$ in $\ltz$ by
\begin{equation}
(Hf)(n) = (L f)(n), \quad n\in\bbZ, \quad
f=\{f(n)\}_{n\in\bbZ}\in\dom(H)=\ltz. \lb{2.4}
\end{equation}
Next, let $g(z,\cdot)$ denote the diagonal Green's function of $H$,
that is,
\begin{equation}
g(z,n)=G(z,n,n), \quad G(z,n,n')=(H-zI)^{-1}(n,n'), \;
z\in\bbC\bs\sigma(H), \; n,n'\in\bbZ.
\end{equation}
Since for each $n\in\bbZ$, $g(\cdot,n)$ is a Herglotz function
(i.e., it maps the open complex upper half-plane analytically to
itself),
\begin{equation}
\xi(\la,n)=\f{1}{\pi}\lim_{\eps\downarrow 0}
\Im[\ln(g(\la+i\eps,n))] \, \text{ for a.e.\ $\la\in\bbR$} \lb{2.5}
\end{equation}
is well-defined for each $n\in\bbZ$. In particular, for all
$n\in\bbZ$,
\begin{equation}
0 \leq \xi(\la,n) \leq 1 \, \text{ for a.e.\ $\la\in\bbR$.}
\end{equation}

In the following we will frequently use the convenient abbreviation
\begin{equation}
h(\la_0+i0)=\lim_{\eps\downarrow 0} h(\la_0 +i\eps), \quad
\la_0\in\bbR, \lb{2.7}
\end{equation}
whenever the limit in \eqref{2.7} is well-defined and hence
\eqref{2.5} can then be written as
$\xi(\la,n)=(1/\pi)\Arg(g(\la+i0,n))$. Moreover, in this section we
will use the convention that whenever the phrase a.e.\ is used
without further qualification, it always refers to Lebesgue measure
on $\bbR$.

Associated with $H$ in $\ltz$, we also introduce the two
half-lattice Jacobi operators $H_{\pm,n_0}$ in
$\lt{[n_0,\pm\infty)\cap\bbZ}$ by
\begin{align}
&H_{\pm,n_0} = P_{\pm,n_0} H
P_{\pm,n_0}|_{\lt{[n_0,\pm\infty)\cap\bbZ}},
\end{align}
where $P_{\pm,n_0}$ are the orthogonal projections onto the
subspaces $\lt{[n_0,\pm\infty)\cap\bbZ}$. By inspection,
$H_{\pm,n_0}$ satisfy Dirchlet boundary conditions at $n_0\mp1$,
that is,
\begin{align}
&(H_{\pm,n_0}f)(n) = (Lf)(n), \quad n\gtreqless n_0, \no
\\
&f\in\dom(H_{\pm,n_0})=\lt{[n_0,\pm\infty)\cap\bbZ}, \quad
f(n_0\mp1)=0.
\end{align}

The half-lattice Weyl--Titchmarsh m-functions associated with
$H_{\pm,n_0}$ are denoted by $m_{\pm}(\cdot,n_0)$ and
$M_{\pm}(\cdot,n_0)$,
\begin{align}
m_\pm(z,n_0) &= (\delta_{n_0}, (H_{\pm,n_0}-zI)^{-1}
\delta_{n_0})_{\lt{[n_0,\pm\infty)\cap\bbZ}}, \quad
z\in\bbC\bs\sigma(H_{\pm,n_0}),   \lb {2.9} \\
M_+(z,n_0) & = -m_+(z,n_0)^{-1}-z+b(n_0), \quad
z\in\bbC\bs\bbR, \lb{2.10}\\
M_-(z,n_0) & = m_-(z,n_0)^{-1}, \quad z\in\bbC\bs\bbR, \lb{2.11}
\end{align}
where $\delta_k = \{\delta_{k,n}\}_{n\in\bbZ}$, $k\in\bbZ$. An
equivalent definition of $M_\pm(\cdot,n_0)$ is
\begin{align}
M_\pm(z,n_0)=-a(n_0)\frac{\psi_\pm(z,n_0+1)}{\psi_\pm(z,n_0)}, \quad
z\in\bbC\bs\bbR, \lb{2.12}
\end{align}
where $\psi_\pm(z,\cdot)$ are the Weyl--Titchmarsh solutions of
$(L-z)\psi_\pm(z,\cdot)=0$ with
$\psi_\pm(z,\cdot)\in\lt{[n_0,\pm\infty)\cap\bbZ}$. Then it follows
that the diagonal Green's function $g(\cdot,n_0)$ is related to the
m-functions $M_\pm(\cdot,n_0)$ via
\begin{align}
g(z,n_0) &= [M_-(z,n_0)-M_+(z,n_0)]^{-1}.   \lb{2.13}
\end{align}
For subsequent purpose we note the universal asymptotic $z$-behavior
of $g(z,n_0)$, valid for all $n_0\in\bbZ$,
\begin{align}
g(z,n_0)\underset{\abs{z}\to\infty}{=}-\frac1z[1+o(1)]. \lb{2.16}
\end{align}

%%%%%%%%%%%%%%%%%%%%%%%%%%%%%%%%%%%%%%%
\begin{definition}  \lb{d2.1}
Let $\cE\subset\bbR$ be a compact set which we may write as
\begin{align}
\cE= [E_0,E_\infty]\Big\bs\bigcup_{j\in J} (E_{2j-1},E_{2j}), \quad
J\subseteq\bbN,
\end{align}
for some $E_0,E_\infty\in\bbR$ and $E_0\leq E_{2j-1}<E_{2j}\leq
E_\infty$, where $(E_{2j-1},E_{2j})\cap
(E_{2j'-1},E_{2j'})=\emptyset$, $j,j'\in J$, $j\neq j'$. Then $\cE$
is called {\it homogeneous} if
\begin{align}
\begin{split}
& \text{there exists an $\eps>0$
such that for all $\la\in\cE$}  \lb{2.18} \\
& \text{and all $0<\delta<\diam(\cE)$, \,
$|\cE\cap(\la-\delta,\la+\delta)|\geq \eps\delta$.}
\end{split}
\end{align}
\end{definition}
%%%%%%%%%%%%%%%%%%%%%%%%%%%%%%%%%%%%%%%

Here $\diam(\cM)$ denotes the diameter of the set $\cM\subset\bbR$.

Next, following \cite{GKT96}, we introduce a special class of
reflectionless Jacobi operators (cf.\ also \cite{Re07} and
\cite[Lemma 8.1]{Te00}).

%%%%%%%%%%%%%%%%%%%%%%%%%%%%%%%%%%%%%%%
\begin{definition}  \lb{d2.2}
Let $\La\subset\bbR$ be of positive Lebesgue measure. Then we call
$H$ reflectionless on $\La$ if for some $n_0\in\bbZ$
\begin{align}
M_+(\la+i0,n_0) = \ol{M_-(\la+i0,n_0)} \, \text{ for a.e. }
\la\in\La. \lb{2.19}
\end{align}
Equivalently $($cf.\ \cite{GS96}$)$, $H$ is called reflectionless on $\La$ if for all
$n\in\bbZ$,
\begin{align}
\xi(\la,n) = 1/2 \, \text{ for a.e. } \la\in\La. \lb{2.20}
\end{align}
\end{definition}
%%%%%%%%%%%%%%%%%%%%%%%%%%%%%%%%%%%%%%%

In the following hypothesis we describe a special class $\cR(\cE)$
of reflectionless Jacobi operators associated with a homogeneous set
$\cE$, that will be our main object of investigation in this
section.

%%%%%%%%%%%%%%%%%%%%%%%%%%%%%%%%%%%%%%%
\begin{hypothesis} \lb{h2.3}
Let $\cE\subset\bbR$ be a compact homogeneous set. Then
$H\in\cR(\cE)$ if
\begin{enumerate}[$(i)$]
\item
$H$ is reflectionless on $\cE$. $($In particular, this implies
$\cE\subseteq\si_{\ess}(H)$.$)$
\item
Either $\cE=\si_{\ess}(H)$ or the set $\si_{\ess}(H)\bs\cE$ is
closed. $($In particular, this implies that there is an open set
$\cO\subseteq\bbR$ such that $\cE\subset\cO$ and
$\ol{\cO}\cap(\si_{\ess}(H)\bs\cE) = \emptyset$.$)$
\item The discrete eigenvalues of $H$ that accumulate to $\cE$
satisfy a Blaschke-type condition, that is,
\begin{align}
\sum_{\la\in\si(H)\cap(\cO\bs\cE)} G_{\cE}(\la,\infty)<\infty,
\lb{2.21}
\end{align}
where $\cO$ is the set defined in $(ii)$ $($hence
$\si(H)\cap(\cO\bs\cE)\subseteq\si_\disc(H)$ is a discrete countable
set$)$ and $G_{\cE}(\cdot,\infty)$ is the potential  theoretic
Green's function for the domain $(\bbC\cup\{\infty\})\bs\cE$ with
logarithmic singularity at infinity $($cf., e.g.,
\cite[Sect.\ 5.2]{Ra95}$)$,
$
G_\cE(z,\infty) \underset{|z|\to\infty}{=} \log |z| - \log (\ca(\cE))
+ \oh(1).
$
\end{enumerate}
\end{hypothesis}
%%%%%%%%%%%%%%%%%%%%%%%%%%%%%%%%%%%%%%%

Here $\ca(\cE)$ denotes the (logarithmic) capacity of $\cE$ (see, e.g.,
\cite[App.\ A]{Si08}).

One particularly interesting situation in which the above hypothesis
is satisfied occurs when $\si(H)=\si_{\ess}(H)$ is a homogeneous set
and $H$ is reflectionless on $\si(H)$. This case has been studied in great
detail by Sodin and Yuditskii \cite{SY97}.

Next, we present the main result of this section. For a Jacobi
operator $H$ in the class $\cR(\cE)$, we will show the absence of
the singular spectrum on the set $\cE$. The proof of this fact
relies on certain techniques developed in harmonic analysis and
potential theory associated with domains
$(\bbC\cup\{\infty\})\bs\cE$ studied by Peherstorfer, Sodin, and
Yuditskii in \cite{PY03} and \cite{SY97}. For completeness, we
provide the necessary result in Theorem \ref{tB.9}.

%%%%%%%%%%%%%%%%%%%%%%%%%%%%%%%%%%%%%%%
\begin{theorem} \lb{t2.7}
Assume Hypothesis \ref{h2.3}, that is, $H\in\cR(\cE)$. Then, the
spectrum of $H$ is purely absolutely continuous on $\cE$,
\begin{equation}
\sigma_{\ac}(H) \supseteq \cE, \quad
\sigma_{\sc}(H)\cap\cE=\sigma_{\pp}(H)\cap\cE=\emptyset. \lb{2.22}
\end{equation}
Moreover, $\sigma(H)$ has uniform multiplicity equal to two on
$\cE$.
\end{theorem}
%%%%%%%%%%%%%%%%%%%%%%%%%%%%%%%%%%%%%%%
\begin{proof}
Fix $n\in\bbZ$. By the asymptotic behavior of the diagonal Green's
function $g(\cdot,n)$ in \eqref{2.16} one concludes that $g(z,n)$ is
a Herglotz function of the type (cf.\ \eqref{B.3} and
\eqref{B.22}--\eqref{B.24})
\begin{equation}
g(z,n)=\int_{\si(H)} \f{d\Om_{0,0}(\la,n)}{\la -z}, \quad
z\in\bbC\bs\si(H).
\end{equation}
Next, we introduce two Herglotz functions $r_j(z,n)$, $j=1,2$, by
\begin{align}
r_1(z,n) &= \int_{\cO} \f{d\Om_{0,0}(\la,n)}{\la-z} =
\sum_{\la\in\cO\bs\cE}\f{\Om_{0,0}(\{\la\},n)}{\la-z} + \int_{\cE}
\f{d\Om_{0,0}(\la,n)}{\la-z}, \quad z\in\bbC\bs\cO,
\lb{2.55} \\
r_2(z,n) &= \int_{\si(H)\bs\cO} \f{d\Om_{0,0}(\la,n)}{\la-z}, \quad
z\in(\bbC\bs\si(H))\cup\cO,
\end{align}
where $\cO$ is the set defined in Hypothesis \ref{h2.3}\,$(ii)$. Then
it is easy to see that
\begin{align}
g(z,n) = r_1(z,n)+r_2(z,n), \quad z\in\bbC\bs\si(H).
\end{align}
Since $H\in\cR(\cE)$, one has $\xi(\la+i0,n)=1/2$ and hence
$\Re[g(\la+i0,n)]=0$ for a.e. $\la\in\cE$. This yields
\begin{align}
\Re[r_1(\la+i0,n)] = -\Re[r_2(\la+i0,n)] \,\text{ for a.e. }
\la\in\cE.
\end{align}
Observing that the function $r_2(\cdot,n)$ is analytic on $(\bbC\bs\si(H))\cup\cO$
and $\cE\subset\cO$, one concludes that $r_2(\cdot,n)$ is bounded on
$\cE$, and hence,
\begin{align}
\Re[r_1(\cdot+i0,n)] = -\Re[r_2(\cdot+i0,n)] \in
L^1\big(\cE;dx\big). \lb{2.59}
\end{align}
Moreover, it follows from Theorem \ref{tB.7} that the set of mass
points of $d\Om_{0,0}$ is a subset of the set of discrete
eigenvalues of $H$, hence \eqref{2.21}, \eqref{2.55}, and
\eqref{2.59} imply that the function $r_1(\cdot,n)$ satisfies the
assumptions of Theorem \ref{tB.9}. Thus, the restriction
$d\Om_{0,0}|_\cE$ of the measure $d\Om_{0,0}$ to the set $\cE$ is
purely absolutely continuous,
\begin{align}
d\Om_{0,0}(\cdot,n)\big|_\cE = d\Om_{0,0,\ac}(\cdot,n)\big|_\cE
=\f{1}{\pi}\Im[r_1(\cdot+i0,n)]d\la\big|_\cE, \quad n\in\bbZ. \lb{2.60}
\end{align}

Finally, utilizing the formulas
\begin{align}
d\Om^{\tr}(\cdot,n) = d\Om_{0,0}(\cdot,n) + d\Om_{1,1}(\cdot,n)
\, \text{ and } \, d\Om_{1,1}(\cdot,n) = d\Om_{0,0}(\cdot,n+1)
\end{align}
one concludes from \eqref{2.60}
\begin{align}
d\Om^{\tr}(\cdot,n)\big|_\cE = d\Om^{\tr}_{\ac}(\cdot,n)\big|_\cE,
\end{align}
for the restriction of the trace measure $d\Om^{\tr}(\cdot,n)$
associated with $H$. By \eqref{B.28} and Theorem \ref{tB.7}\,$(i)$
this completes the proof of \eqref{2.22}.

Finally, equations \eqref{2.13} and \eqref{2.19} imply
\begin{equation}
-1/g(\la+i0,n)=\pm 2i \, \Im[M_\pm(\la+i0,n)] \, \text{ for a.e.\
$\la\in\cE$.}    \lb{2.30}
\end{equation}
Thus, combining \eqref{2.19}, \eqref{2.30}, and \eqref{B.38} then
yields that the absolutely continuous spectrum of $H$ has uniform
spectral multiplicity two on $\cE$ since
\begin{equation}
\text{for a.e.\ $\la\in\cE$, } \; 0<\pm\Im[M_\pm(\la+i0,n)]<\infty.
\end{equation}
\end{proof}
%%%%%%%%%%%%%%%%%%%%%%%%%%%%%%%%%%%%%%%

For reflectionless measures with singular components we refer to the recent
preprint \cite{NVY08} (see also \cite{GY06}).

%%%%%%%%%%%%%%%%%%%%%%%%%%%%%%%%%%%%%%%
%%%%%%%%%%%%%%%%%%%%%%%%%%%%%%%%%%%%%%%
\section{Reflectionless CMV Operators} \lb{s3}
%%%%%%%%%%%%%%%%%%%%%%%%%%%%%%%%%%%%%%%
%%%%%%%%%%%%%%%%%%%%%%%%%%%%%%%%%%%%%%%

In this section we investigate spectral properties of unitary CMV
operators reflectionless on compact homogeneous subsets of the unit
circle.

We start with some general considerations of unitary CMV operators.
Let $\{\al_n\}_{n\in\bbZ}$ be a complex-valued sequence of
Verblunsky coefficients satisfying
\begin{equation}
\al_n\in\bbD = \{z\in\bbC \st |z|<1\}, \quad n\in\bbZ,
\end{equation}
and denote by $\{\rho_n\}_{n\in\bbZ}$ an auxiliary real-valued
sequence defined by
\begin{align}
\rho_n = \big[1-\abs{\al_n}^2\big]^{1/2}, \quad n\in\bbZ.
\end{align}
Then we introduce the associated unitary CMV operator $U$ in $\ltz$
by its matrix representation in the standard basis of $\ltz$,
\begin{align}
U = \begin{pmatrix} \ddots &&\hspace*{-8mm}\ddots
&\hspace*{-10mm}\ddots &\hspace*{-12mm}\ddots &\hspace*{-14mm}\ddots
&&& \raisebox{-3mm}[0mm][0mm]{\hspace*{-6mm}{\Huge $0$}}
\\
&0& -\al_{0}\rho_{-1} & -\ol{\al_{-1}}\al_{0} & -\al_{1}\rho_{0} &
\rho_{0}\rho_{1}
\\
&& \rho_{-1}\rho_{0} &\ol{\al_{-1}}\rho_{0} & -\ol{\al_{0}}\al_{1} &
\ol{\al_{0}}\rho_{1} & 0
\\
&&&0& -\al_{2}\rho_{1} & -\ol{\al_{1}}\al_{2} & -\al_{3}\rho_{2} &
\rho_{2}\rho_{3}
\\
&&\raisebox{-4mm}[0mm][0mm]{\hspace*{-6mm}{\Huge $0$}} &&
\rho_{1}\rho_{2} & \ol{\al_{1}}\rho_{2} & -\ol{\al_{2}}\al_{3} &
\ol{\al_{2}}\rho_{3}&0
\\
&&&&&\hspace*{-14mm}\ddots &\hspace*{-14mm}\ddots
&\hspace*{-14mm}\ddots &\hspace*{-8mm}\ddots &\ddots
\end{pmatrix}. \lb{3.3}
\end{align}
Here terms of the form $-\ol{\alpha_n} \alpha_{n+1}$ represent the
diagonal $(n,n)$-entries, $n\in\bbZ$, in the infinite matrix
\eqref{3.3}. Equivalently, one can define $U$ by (cf.\ \eqref{2.3})
\begin{align}
U &= \rho^- \rho \, \deven \, S^{--} + (\ol{\alpha^-}\rho \, \deven -
\alpha^+\rho \, \dodd) S^-
- \ol\alpha\alpha^+   \no \\
& \quad + (\ol\alpha \rho^+ \, \deven - \alpha^{++} \rho^+ \, \dodd)
S^+ + \rho^+ \rho^{++} \, \dodd \, S^{++},    \lb{3.4}
\end{align}
where $\deven$ and $\dodd$ denote the characteristic functions of the
even and odd integers,
\begin{equation}
\deven = \chi_{_{2\bbZ}}, \quad \dodd = 1 - \deven = \chi_{_{2\bbZ
+1}}.
\end{equation}

Moreover, let $M_{1,1}(z,n)$ denote the diagonal element of the
Cayley transform of $U$, that is,
\begin{align}
M_{1,1}(z,n)=((U+zI)(U-zI)^{-1})(n,n) = \oint_\dD
d\Om_{1,1}(\ze,n)\,\f{\ze+z}{\ze-z},& \no
\\
z\in\bbC\bs\sigma(U), \; n\in\bbZ,&
\end{align}
where $d\Om_{1,1}(\cdot,n)$, $n\in\bbZ$, are scalar-valued
probability measures on $\dD$ (cf.\ \cite[Section 3]{GZ06} for more
details). Since $M_{1,1}(\cdot,n)$ is a Caratheodory function (i.e.,
it maps the open unit disk analytically to the complex right
half-plane),
\begin{equation}
\Xi_{1,1}(\ze,n)=\f{1}{\pi}\lim_{r\uparrow 1}
\Im[\ln(M_{1,1}(r\ze,n))] \, \text{ for a.e.\ $\ze\in\dD$} \lb{3.5}
\end{equation}
is well-defined for each $n\in\bbZ$. In particular, for all
$n\in\bbZ$,
\begin{equation}
-1/2 \leq \Xi_{1,1}(\ze,n) \leq 1/2 \, \text{ for a.e.\ $\ze\in\dD$}
\end{equation}
(cf.\ \cite[Section 2]{GZ06a} for more details).

In the following we will frequently use the convenient abbreviation
\begin{equation}
h(\ze)=\lim_{r\uparrow 1} h(r\ze), \quad \ze\in\dD, \lb{3.7}
\end{equation}
whenever the limit in \eqref{3.7} is well-defined and hence
\eqref{3.5} can then be written as
$\Xi_{1,1}(\ze,n)=(1/\pi)\Arg(M_{1,1}(\ze,n))$. Moreover, in this
section we will use the convention that whenever the phrase a.e.\ is
used without further qualification, it always refers to Lebesgue
measure on $\dD$.

Associated with $U$ in $\ltz$, we also introduce the two half-lattice
CMV operators $U_{\pm,n_0}$ in $\lt{[n_0,\pm\infty)\cap\bbZ}$ by
setting $\al_{n_0}=1$ which splits the operator $U$ into a direct sum
of two half-lattice operators $U_{-,n_0-1}$ and $U_{+,n_0}$, that is,
\begin{align}
U=U_{-,n_0-1} \oplus U_{+,n_0} \, \text{ in } \,
\lt{(-\infty,n_0-1]\cap\bbZ} \oplus \lt{[n_0,\infty)\cap\bbZ}.
\end{align}
The half-lattice Weyl--Titchmarsh m-functions associated with
$U_{\pm,n_0}$ are denoted by $m_{\pm}(\cdot,n_0)$ and
$M_{\pm}(\cdot,n_0)$,
\begin{align}
m_\pm(z,n_0) &= ((U_{\pm,n_0}+zI)(U_{\pm,n_0}-zI)^{-1})(n_0,n_0),
\quad
z\in\bbC\bs\sigma(U_{\pm,n_0}), \\
M_+(z,n_0) & = m_+(z,n_0), \quad
z\in\bbC\bs\dD, \lb{3.10}\\
M_-(z,n_0) &= \frac{\Re(1+\al_{n_0}) +
i\Im(1-\al_{n_0})m_-(z,n_0-1)}{i\Im(1+\al_{n_0}) +
\Re(1-\al_{n_0})m_-(z,n_0-1)}, \quad z\in\bbC\bs\dD. \lb{3.11}
\end{align}
Then it follows that $m_\pm(\cdot,n_0)$ and $\pm M_\pm(\cdot,n_0)$
are Caratheodory functions (cf.\ \cite[Section 2]{GZ06}). Moreover,
the function $M_{1,1}(\cdot,n_0)$ is related to the m-functions
$M_\pm(\cdot,n_0)$ by (cf.\ \cite[Lemma 3.2]{GZ06})
\begin{align}
M_{1,1}(z,n_0) &=
\frac{1-M_+(z,n_0)M_-(z,n_0)}{M_+(z,n_0)-M_-(z,n_0)}. \lb{3.13}
\end{align}

%%%%%%%%%%%%%%%%%%%%%%%%%%%%%%%%%%%%%%%
\begin{definition}  \lb{d3.1}
Let $\cE\subseteq\dD$ be a compact set which we may write as
\begin{align}
\cE=\dD\big\bs\bigcup_{j\in J} \Arc(e^{i\te_{2j-1}},e^{i\te_{2j}}),
\quad J\subseteq\bbN,
\end{align}
where $\Arc(e^{i\te_{2j-1}},e^{i\te_{2j}}) =\big\{e^{i\te}\in\dD \st
\te_{2j-1}<\theta<\te_{2j}\big\}$, $\theta_{2j-1}\in[0,2\pi)$,
$\te_{2j-1}<\te_{2j}\leq\te_{2j-1}+2\pi$,
$\Arc(e^{i\te_{2j-1}},e^{i\te_{2j}})\cap
\Arc(e^{i\te_{2j'-1}},e^{i\te_{2j'}})=\emptyset$, $j,j'\in J$,
$j\neq j'$. Then $\cE$ is called {\it homogeneous} if
\begin{align}
\begin{split}
& \text{there exists an $\eps>0$
such that for all $e^{i\te}\in\cE$}  \lb{3.18} \\
& \text{and all $\de>0$, \,
$|\cE\cap\Arc(e^{i(\te-\de)},e^{i(\te+\de)})|\geq \eps\delta$.}
\end{split}
\end{align}
\end{definition}
%%%%%%%%%%%%%%%%%%%%%%%%%%%%%%%%%%%%%%%

Next, we introduce a special class of reflectionless CMV operators
(cf.\ \cite{GZ06a} for a similar definition).

%%%%%%%%%%%%%%%%%%%%%%%%%%%%%%%%%%%%%%%
\begin{definition}  \lb{d3.2}
Let $\La\subseteq\dD$ be of positive Lebesgue measure. Then we call
$U$ reflectionless on $\La$ if for some $($equivalently, for all\,$)$
$n_0\in\bbZ$
\begin{align}
M_+(\ze,n_0) = -\ol{M_-(\ze,n_0)} \, \text{ for a.e.\ } \ze\in\La.
\lb{3.19}
\end{align}
\end{definition}
%%%%%%%%%%%%%%%%%%%%%%%%%%%%%%%%%%%%%%%

We note that if $U$ is reflectionless on $\La$, then by \eqref{3.5},
\eqref{3.13}, and \eqref{3.19}, one has for all $n\in\bbZ$,
\begin{align}
\Xi_{1,1}(\ze,n) = 0 \, \text{ for a.e. } \ze\in\La. \lb{3.20}
\end{align}

In the following hypothesis we introduce a special class $\cR(\cE)$
of reflectionless CMV operators associated with a homogeneous set
$\cE$, that will be the main object of investigation in this section.

%%%%%%%%%%%%%%%%%%%%%%%%%%%%%%%%%%%%%%%
\begin{hypothesis} \lb{h3.3}
Let $\cE\subseteq\dD$ be a compact homogeneous set. Then
$U\in\cR(\cE)$ if
\begin{enumerate}[$(i)$]
\item
$U$ is reflectionless on $\cE$. $($In particular, this implies
$\cE\subseteq\si_{\ess}(U)$.$)$
\item
Either $\cE=\si_{\ess}(U)$ or the set $\si_{\ess}(U)\bs\cE$ is
closed. $($In particular, this implies that there is an open set
$\cO\subseteq\dD$ such that $\cE\subseteq\cO$ and
$\ol{\cO}\cap(\si_{\ess}(U)\bs\cE) = \emptyset$.$)$
\item If $\si_\ess(U)\neq\dD$ then the discrete eigenvalues of\, $U$
that accumulate to $\cE$ satisfy a Blaschke-type condition, that is,
\begin{align}
\sum_{\ze\in\si(U)\cap(\cO\bs\cE)} G_{\cE}(\ze,\ze_0)<\infty, \lb{3.21}
\end{align}
where $\cO$ is the set defined in $(ii)$ $($hence
$\si(U)\cap(\cO\bs\cE)\subseteq\si_\disc(U)$ is a discrete countable
set$)$, $\ze_0\in\dD\bs\si(U)$ is some fixed point, and
$G_{\cE}(\cdot,\ze_0)$ is the potential  theoretic Green's function
for the domain $(\bbC\cup\{\infty\})\bs\cE$ with logarithmic
singularity at $\ze_0$ $($cf., e.g., \cite[Sect.\ 4.4]{Ra95}$)$,
$
G_\cE(z,\ze_0) \underset{z\to \ze_0}{=} \log |z-\ze_0|^{-1} +
\Oh(1).
$
\end{enumerate}
\end{hypothesis}
%%%%%%%%%%%%%%%%%%%%%%%%%%%%%%%%%%%%%%%
One particularly interesting situation in which the above hypothesis
is satisfied occurs when $\si(U)=\si_{\ess}(U)$ is a homogeneous set
and $U$ is reflectionless on $\si(U)$. This case has first been studied by
Peherstorfer and Yudiskii \cite{PY06}.

Next, we turn to the principal result of this section. For a CMV operator
$U$ in the class $\cR(\cE)$, we will show that $U$ has purely
absolutely continuous spectrum on $\cE$, that is, we intend to prove
that
\begin{equation}
\sigma_{\ac}(U)\supseteq\cE, \quad
\sigma_{\sc}(U)\cap\cE=\sigma_{\pp}(U)\cap\cE=\emptyset.
\end{equation}

We start with an elementary lemma which permits one to apply Theorem
\ref{tB.9} to Caratheodory functions.

%%%%%%%%%%%%%%%%%%%%%%%%%%%%%%%%%%%%%%%
\begin{lemma}  \lb{l3.4}
Let $f$ be a Caratheodory function with representation
\begin{align}
\begin{split}
& f(w)=ic+ \oint_{\dD} d\om(\ze) \, \f{\ze+w}{\ze-w}, \quad
w\in\bbD,
\\
& c=\Im(f(0)), \;\; \oint_{\dD} d\om(\ze) = \Re(f(0)) < \infty, \;\;
\supp\,(d\om)\neq\dD,
\end{split}
\end{align}
where $d\om$ denotes a nonnegative measure on $\dD$. Consider the
change of variables
\begin{align}
& w \mapsto z=-i\f{w+w_0}{w-w_0}, \;\; w=w_0\f{z-i}{z+i},
\quad z\in\bbC\cup\{\infty\}, \\
& \text{for some fixed $w_0\in\dD\bs\supp\,(d\om)$}. \no
\end{align}
Then, the function $r(z)=if(w(z))$ is a Herglotz function with the
representation
\begin{align}
& r(z)=if(w(z))=d+\int_{\bbR} \f{d\mu(\la)}{\la-z}, \quad
z\in\bbC_+, \lb{3.50}
\\
& d=if(w_0), \;\; d\mu(\la)=(1+\la^2)\,d\om(w(\la)), \quad
\la\in\bbR, \lb{3.51}
\\
& \supp\,(d\mu) \subseteq \big[-1-2\dist(w_0,\supp\,(d\om))^{-1},
1+2\dist(w_0,\supp\,(d\om))^{-1}\big].
\end{align}
In particular, $d\mu$ is purely absolutely continuous on $\La$ if
and only if $d\om$ is purely absolutely continuous on $w(\La)$, and
\begin{equation}
\text{if\, $d\om(e^{i\te})\big|_{w(\La)}=\om'(e^{i\te})d\te\big|_{w(\La)}$,\, then\,
$d\mu(\la)\big|_{\Lambda}=2\om'(w(\la))\,d\la\big|_{\Lambda}$.}
\end{equation}
\end{lemma}
%%%%%%%%%%%%%%%%%%%%%%%%%%%%%%%%%%%%%%
\begin{proof}
This is a straightforward computation. We note that
\begin{equation}
\oint_{\dD}d\om(\ze)<\infty \, \text{ is equivalent to }
\int_{\bbR}\f{d\mu(\la)}{1+\la^2}<\infty.
\end{equation}
\end{proof}
%%%%%%%%%%%%%%%%%%%%%%%%%%%%%%%%%%%%%%

The principal result of this section then reads as follows.

%%%%%%%%%%%%%%%%%%%%%%%%%%%%%%%%%%%%%%%
\begin{theorem} \lb{t3.5}
Assume Hypothesis \ref{h3.3}, that is, $U\in\cR(\cE)$. Then, the
spectrum of $U$ is purely absolutely continuous on $\cE$,
\begin{equation}
\sigma_{\ac}(U) \supseteq \cE, \quad
\sigma_{\sc}(U)\cap\cE=\sigma_{\pp}(U)\cap\cE=\emptyset. \lb{3.53}
\end{equation}
Moreover, $\sigma(U)$ has uniform multiplicity equal to two on $\cE$.
\end{theorem}
%%%%%%%%%%%%%%%%%%%%%%%%%%%%%%%%%%%%%%%
\begin{proof}
We consider two cases. First, suppose that
$\si_\ess(U)\neq\dD$. Then using Lemma \ref{l3.4}, we introduce the
Herglotz function $r(\cdot,n)$, $n\in\bbZ$, by
\begin{align}
r(z,n) = iM_{1,1}(w(z),n) =
iM_{1,1}(\ze_0,n)+\int_\bbR\f{d\mu(\la,n)}{\la-z}, \quad z\in\bbC_+,
\end{align}
where $\ze_0$ is defined in Hypothesis \ref{h3.3}\,$(iii)$ and
$w(z)=\ze_0\f{z-i}{z+i}$. Abbreviating by $\wti\cE$ and $\wti\cO$
the preimages of the sets $\cE$ and $\cO$ under the bijective map
$w$,
\begin{align}
\wti\cE=w^{-1}(\cE), \quad \wti\cO=w^{-1}(\cO),
\end{align}
we also introduce functions $r_j(z,n)$, $j=1,2$, by
\begin{align}
r_1(z,n) &= \int_{\wti\cO} \f{d\mu(\la,n)}{\la-z} = \sum_{\la\in
\wti\cO\bs\wti\cE}\f{\mu(\{\la\},n)}{\la-z} + \int_{\wti\cE}
\f{d\mu(\la,n)}{\la-z}, \quad z\in\bbC\bs\wti\cO,
\lb{3.55} \\
r_2(z,n) &= \int_{\bbR\bs\wti\cO} \f{d\mu(\la,n)}{\la-z}, \quad
z\in(\bbC\bs\bbR)\cup\wti\cO.
\end{align}
Then
\begin{align}
r(z,n) = r_1(z,n)+r_2(z,n), \quad z\in\bbC\bs\bbR, \; n\in\bbZ.
\end{align}

Since $U\in\cR(\cE)$, one has $\Xi_{1,1}(\ze,n)=0$ and hence
$\Im[M_{1,1}(\ze,n)]=0$ for all $n\in\bbZ$ and a.e. $\ze\in\cE$.
This yields for each $n\in\bbZ$, $\Re[r(\la+i0,n)] = 0$ for a.e.
$\la\in\wti\cE$, and hence,
\begin{align}
\Re[r_1(\la+i0,n)] = -\Re[r_2(\la+i0,n)] \,\text{ for a.e. }
\la\in\wti\cE.
\end{align}
Observing that the function $r_2(\cdot,n)$ is analytic on $(\bbC\bs\bbR)\cup\wti\cO$
and $\wti\cE\subset\wti\cO$, one concludes that $r_2(\cdot,n)$ is
bounded on $\wti\cE$, and hence,
\begin{align}
\Re[r_1(\cdot+i0,n)] = -\Re[r_2(\cdot+i0,n)] \in
L^1\big(\wti\cE;dx\big), \quad n\in\bbZ. \lb{3.59}
\end{align}
Moreover, it follows from \cite[Proposition 5.1]{Fi83} that
\eqref{3.21} is equivalent to
\begin{align}
\sum_{\la\in\wti\cO\bs\wti\cE} G_{\wti\cE}(\la,\infty)<\infty,
\lb{3.59a}
\end{align}
and from \cite[Corollary 3.5]{GZ06} that the set of discrete mass
points of $d\mu(\cdot,n)$ is a subset of $\wti\cO\bs\wti\cE$, hence
\eqref{3.55}, \eqref{3.59}, and \eqref{3.59a} imply that the
function $r_1(\cdot,n)$ satisfies the assumptions of Theorem
\ref{tB.9} for each $n\in\bbZ$. Thus, the restriction
$d\mu(\cdot,n)|_{\wti\cE}$ of the measure $d\mu(\cdot,n)$ to the set
$\wti\cE$ is purely absolutely continuous,
\begin{align}
& d\mu(\cdot,n)\big|_{\wti\cE} = d\mu_\ac(\cdot,n)\big|_{\wti\cE}
=\f{1}{\pi}\Im[r_1(\cdot+i0,n)]d\la\big|_{\wti\cE}, \quad n\in\bbZ,
\end{align}
and hence, it follows from Lemma \ref{l3.4} that
\begin{align}
d\Om_{1,1}(\cdot,n)\big|_\cE = d\Om_{1,1,\ac}(\cdot,n)\big|_\cE,
\quad n\in\bbZ. \lb{3.60}
\end{align}
By Theorem \ref{tC.7}, and in particular \eqref{C.36a}, this proves
\eqref{3.53} in the case $\si_\ess(U)\neq\dD$.

Next, suppose $\si_\ess(U)=\dD$. In this case it follows from a
special case of the Borg-type theorem proven in \cite[Theorem 5.1]{GZ06a}
(cf.\ also \cite[Sect.\ 11.14]{Si05} for a more
restrictive version of this theorem) that the Verblunsky coefficients
$\al_n=0$ for all $n\in\bbZ$. Hence, \eqref{3.4} implies that $U$ is
unitarily equivalent to a direct sum of two shift operators in $\ltz$
($U$ shifts odd entries to the left and even entries to the right).
Thus $U$ has purely absolutely continuous spectrum on $\dD$, which
proves \eqref{3.53} in the case $\si_\ess(U)=\dD$.

Finally, equations \eqref{3.13} and \eqref{3.19} imply
\begin{equation}
\frac{1}{M_{1,1}(\ze,n)}=\frac{\pm
2\,\Re[M_\pm(\ze,n)]}{1+|M_\pm(\ze,n)|^2} \, \text{ for a.e.\
$\ze\in\cE$.} \lb{3.30}
\end{equation}
Combining \eqref{3.19}, \eqref{3.30}, and \eqref{C.38} then yields
that the absolutely continuous spectrum of $U$ has uniform spectral
multiplicity two on $\cE$ since
\begin{equation}
\text{for a.e.\ $\ze\in\cE$, } \; 0<\pm\Re[M_\pm(\ze,n)]<\infty.
\end{equation}
\end{proof}
%%%%%%%%%%%%%%%%%%%%%%%%%%%%%%%%%%%%%%%

%%%%%%%%%%%%%%%%%%%%%%%%%%%%%%%%%%%%%%%%
%%%%%%%%%%%%%%%%%%%%%%%%%%%%%%%%%%%%%%%%
\section{Reflectionless Schr\"odinger Operators}    \label{s4}
%%%%%%%%%%%%%%%%%%%%%%%%%%%%%%%%%%%%%%%%
%%%%%%%%%%%%%%%%%%%%%%%%%%%%%%%%%%%%%%%%

In this section we discuss spectral properties of self-adjoint Schr\"odinger
operators reflectionless on a homogeneous subsets of the real line
bounded from below.

We start with some general considerations of one-dimensional
Schr\"odinger operators. Let
\begin{equation}
V\in L^\infty(\bbR;dx), \quad V \, \text{ real-valued,}
\end{equation}
and consider the differential expression
\begin{equation}
L=-d^2/dx^2+V(x), \quad x\in\bbR.
\end{equation}
We denote by $H$ the corresponding self-adjoint realization of $L$
in $L^2(\bbR;dx)$ given by
\begin{equation}
Hf=L f, \quad f\in\dom(H)=H^2(\bbR),   \lb{4.4}
\end{equation}
with $H^2(\bbR)$ the usual Sobolev space. Let $g(z,\cdot)$ denote
the diagonal Green's function of $H$, that is,
\begin{equation}
g(z,x)=G(z,x,x), \quad G(z,x,x')=(H-zI)^{-1}(x,x'), \quad
z\in\bbC\backslash\sigma(H), \; x,x'\in\bbR.
\end{equation}
Since for each $x\in\bbR$, $g(\cdot,x)$ is a Herglotz function,
\begin{equation}
\xi(\lambda,x)=\f{1}{\pi}\lim_{\varepsilon\downarrow 0}
\Im[\ln(g(\lambda+i\varepsilon,x))] \, \text{ for a.e.\
$\lambda\in\bbR$}  \lb{4.5}
\end{equation}
is well-defined for each $x\in\bbR$. In particular, for all
$x\in\bbR$,
\begin{equation}
0 \leq \xi(\lambda,x) \leq 1 \, \text{ for a.e.\ $\lambda\in\bbR$.}
\end{equation}

In the following we will frequently use the convenient abbreviation
\begin{equation}
h(\lambda_0+i0)=\lim_{\varepsilon\downarrow 0} h(\lambda_0
+i\varepsilon), \quad \lambda_0\in\bbR, \lb{4.7}
\end{equation}
whenever the limit in \eqref{4.7} is well-defined and hence
\eqref{4.5} can then be written as
$\xi(\lambda,x)=(1/\pi)\Arg(g(\lambda+i0,x))$. Moreover, in this
section we will use the convention that whenever the phrase a.e.\ is
used without further qualification, it always refers to Lebesgue
measure on $\bbR$.

Associated with $H$ in $L^2(\bbR;dx)$ we also introduce the two
half-line Schr\"odinger operators $H_{\pm,x_0}$ in
$L^2([x_0,\pm\infty);dx)$ with Dirchlet boundary conditions at the
finite endpoint $x_0\in\bbR$,
\begin{align}
& H_{\pm,x_0} f= Lf, \no \\
& f\in\dom(H_{\pm,x_0})=\big\{g\in L^2([x_0,\pm\infty);dx)\,|\, g,
g'\in
AC([x_0, x_0\pm R]) \, \text{for all $R>0$;} \no \\
& \hspace*{4.5cm} \lim_{\varepsilon\downarrow 0}
g(x_0\pm\varepsilon)=0; \, Lg\in L^2([x_0,\pm\infty);dx)\big\}.
\end{align}
Denoting by $\psi_{\pm}(z,\cdot)$ the Weyl--Titchmarsh solutions of
$(L - z)\psi(z,\dott) = 0$, satisfying
\begin{equation}
\psi_{\pm}(z,\cdot)\in L^2([x_0,\pm\infty);dx),
\end{equation}
the half-line Weyl--Titchmarsh functions associated with
$H_{\pm,x_0}$ are given by
\begin{equation}
m_{\pm}(z,x_0)=\f{\psi_{\pm}'(z,x_0)}{\psi_{\pm}(z,x_0)}, \quad
z\in\bbC\backslash\sigma(H_{\pm,x_0}). \lb{4.10}
\end{equation}
Then the diagonal Green's function of $H$ satisfies
\begin{equation}
g(z,x_0)=[m_{-}(z,x_0) - m_{+}(z,x_0)]^{-1}.    \lb{4.11}
\end{equation}
For subsequent purpose we also introduce two Herglotz functions,
\begin{align}
M_{0,0}(z,x_0) = \frac{1}{m_{-}(z,x_0) - m_{+}(z,x_0)} = \int_\bbR
\frac{d\Om_{0,0}(\la,x_0)\,d\la}{\la-z}, \lb{4.11a} \\
M_{1,1}(z,x_0) = \frac{m_{-}(z,x_0)m_{+}(z,x_0)}{m_{-}(z,x_0) -
m_{+}(z,x_0)} = \int_\bbR \frac{d\Om_{1,1}(\la,x_0)\,d\la}{\la-z}.
\lb{4.11b}
\end{align}

%%%%%%%%%%%%%%%%%%%%%%%%%%%%%%%%%%%%%%%%
\begin{definition}  \lb{d5.1}
Let $\cE\subset\bbR$ be a closed set bounded from below which we may write
as
\begin{equation}
\cE= [E_0,\infty)\big\backslash\bigcup_{j\in J} (a_j,b_j), \quad
J\subseteq\bbN, \lb{4.12}
\end{equation}
for some $E_0\in\bbR$ and $a_j<b_j$, where $(a_j,b_j)\cap
(a_{j'},b_{j'})=\emptyset$, $j,j'\in J$, $j\neq j'$. Then $\cE$ is called
{\it homogeneous} if
\begin{align}
\begin{split}
& \text{there exists an $\varepsilon>0$
such that for all $\lambda\in\cE$} \\
& \text{and all $\delta>0$, \,
$|\cE\cap(\lambda-\delta,\lambda+\delta)|\geq \varepsilon\delta$.}
\end{split}
\end{align}
\end{definition}
%%%%%%%%%%%%%%%%%%%%%%%%%%%%%%%%%%%%%%%%

Next, we introduce a special class of reflectionless Schr\"odinger
operators.

%%%%%%%%%%%%%%%%%%%%%%%%%%%%%%%%%%%%%%%
\begin{definition}  \lb{d4.2}
Let $\La\subset\bbR$ be of positive Lebesgue measure. Then we call
$H$ reflectionless on $\La$ if for some $x_0\in\bbR$
\begin{equation}
m_+(\lambda+i0,x_0)= \ol{m_-(\lambda+i0,x_0)} \, \text{ for a.e.\
$\lambda\in\La$.}  \lb{4.19}
\end{equation}
Equivalently $($cf.\ \cite{GS96}, \cite{GY06}$)$, $H$ is called reflectionless on $\La$ if for each $x\in\bbR$,
\begin{equation}
\xi(\lambda,x)=1/2 \, \text{ for a.e.\ $\lambda\in\La$.} \lb{4.20}
\end{equation}
\end{definition}
%%%%%%%%%%%%%%%%%%%%%%%%%%%%%%%%%%%%%%%

In the following hypothesis we describe a special class $\cR(\cE)$
of reflectionless Schr\"odinger operators associated with a
homogeneous set $\cE$, that will be our main object of investigation
in this section.

%%%%%%%%%%%%%%%%%%%%%%%%%%%%%%%%%%%%%%%
\begin{hypothesis} \lb{h4.3}
Let $\cE\subset\bbR$ be a homogeneous set. Then $H\in\cR(\cE)$ if
\begin{enumerate}[$(i)$]
\item
$H$ is reflectionless on $\cE$. $($In particular, this implies
$\cE\subseteq\si_{\ess}(H)$.$)$
\item
Either $\cE=\si_{\ess}(H)$ or the set $\si_{\ess}(H)\bs\cE$ is
closed. $($In particular, this implies that there is an open set
$\cO\subseteq\bbR$ such that $\cE\subset\cO$ and
$\ol{\cO}\cap(\si_{\ess}(H)\bs\cE) = \emptyset$.$)$
\item The discrete eigenvalues of $H$ that accumulate to $\cE$
satisfy a Blaschke-type condition, that is,
\begin{align}
\sum_{\la\in\si(H)\cap(\cO\bs\cE)} G_{\cE}(\la,\la_0)<\infty,
\lb{4.21}
\end{align}
where $\cO$ is the set defined in $(ii)$ $($hence
$\si(H)\cap(\cO\bs\cE)\subseteq\si_\disc(H)$ is a discrete countable
set$)$, $\la_0\in\bbR\bs\si(H)$ is some fixed point, and
$G_{\cE}(\cdot,\infty)$ is the potential  theoretic Green's function
for the domain $(\bbC\cup\{\infty\})\bs\cE$ with logarithmic
singularity at $\la_0$ $($cf., e.g., \cite[Sect.\ 4.4]{Ra95}$)$, $
G_\cE(z,\la_0) \underset{z\to\la_0}{=} \log |z-\la_0|^{-1} + \Oh(1).
$
\end{enumerate}
\end{hypothesis}
%%%%%%%%%%%%%%%%%%%%%%%%%%%%%%%%%%%%%%%

One particularly interesting situation in which the above hypothesis
is satisfied occurs when $\si(H)=\si_{\ess}(H)$ is a homogeneous set
and $H$ is reflectionless on $\si(H)$. This case has been studied in
great detail by Sodin and Yuditskii \cite{SY95a}--\cite{SY96}, and
more recently, in \cite{GY06}.

Next, we turn to the principal result of this section. For a
Schr\"odinger operator $H$ in the class $\cR(\cE)$, we will show
that $H$ has purely absolutely continuous spectrum on $\cE$, that
is, we intend to prove that
\begin{equation}
\sigma_{\ac}(H)\supseteq\cE, \quad
\sigma_{\sc}(H)\cap\cE=\sigma_{\pp}(H)\cap\cE=\emptyset.
\end{equation}

We start with an elementary lemma which will permit us to reduce the
discussion of unbounded homogeneous sets $\cE$ (typical for
Schr\"odinger operators) to the case of compact homogeneous sets
$\wti\cE$ (typical for Jacobi operators).

%%%%%%%%%%%%%%%%%%%%%%%%%%%%%%%%%%%%%%%
\begin{lemma}  \lb{l4.4}
Let $m$ be a Herglotz function with representation
\begin{align}
\begin{split}
& m(z)=c+\int_{\bbR} d\omega(\lambda)
\bigg(\f{1}{\lambda-z}-\f{\lambda}{1+\lambda^2}\bigg), \quad
z\in\bbC_+,
\\
& c=\Re[m(i)], \quad \int_{\bbR} \f{d\omega(\lambda)} {1+\lambda^2}
<\infty, \quad \bbR\backslash \supp\,(d\omega)\neq\emptyset.
\end{split}
\end{align}
Consider the change of variables
\begin{align}
\begin{split}
& z \mapsto \zeta=(\lambda_0-z)^{-1}, \;\;
z=\lambda_0-\zeta^{-1}, \quad z\in\bbC\cup\{\infty\}, \\
& \text{for some fixed $\lambda_0\in \bbR\backslash\supp \,
(d\omega)$}.
\end{split}
\end{align}
Then, the function $r(\ze) = m(z(\ze))$ is a Herglotz function with
the representation
\begin{align}
& r(\zeta)=m(z(\zeta)) =d +\int_{\bbR} \f{d\mu(\eta)}{\eta-\zeta},
\quad
\zeta\in\bbC_+, \lb{4.50} \\
& d\mu(\eta)
=x^2d\omega\big(\la_0-\eta^{-1}\big)\big|_{\supp\,(d\mu)},
\lb{4.51} \\
& \supp\, (d\mu) \subseteq
\big[-\dist\,(\lambda_0,\supp\,(d\omega))^{-1},
\dist\,(\lambda_0,\supp\,(d\omega))^{-1}\big], \\
& d=c+\int_{\supp\,(d\mu)} d\mu(\eta)\,
\f{\lambda_0-\big(1+\lambda_0^2\big)\eta} {1-2\lambda_0
\eta+\big(1+\lambda_0^2\big)\eta^2}. \lb{4.52}
\end{align}
In particular, $d\mu$ is purely absolutely continuous on $\La$ if
and only if $d\omega$ is purely absolutely continuous on $z(\La)$,
and
\begin{equation}
\text{if }\,
d\omega(\lambda)|_{z(\La)}=\omega'(\lambda)d\lambda\big|_{z(\La)},
\text{ then }\, d\mu(\eta)|_{\La} =
\omega'\big(\lambda_0-\eta^{-1}\big)d\eta\big|_{\La}.
\end{equation}
\end{lemma}
%%%%%%%%%%%%%%%%%%%%%%%%%%%%%%%%%%%%%%
\begin{proof}
This is a straightforward computation. We note that
\begin{equation}
\int_{\bbR} \f{d\omega(\lambda)}{1+\lambda^2} <\infty \, \text{ is
equivalent to } \int_{\bbR} d\mu(\eta)<\infty.
\end{equation}
\end{proof}
%%%%%%%%%%%%%%%%%%%%%%%%%%%%%%%%%%%%%%

The principal result of this section then reads as follows.

%%%%%%%%%%%%%%%%%%%%%%%%%%%%%%%%%%%%%%%
\begin{theorem} \lb{t4.5}
Assume Hypothesis \ref{h4.3}, that is, $H\in\cR(\cE)$. Then, the
spectrum of $H$ is purely absolutely continuous on $\cE$,
\begin{equation}
\sigma_{\ac}(H) \supseteq \cE, \quad
\sigma_{\sc}(H)\cap\cE=\sigma_{\pp}(H)\cap\cE=\emptyset. \lb{4.53}
\end{equation}
Moreover, $\sigma(H)$ has uniform multiplicity equal to two on
$\cE$.
\end{theorem}
%%%%%%%%%%%%%%%%%%%%%%%%%%%%%%%%%%%%%%%
\begin{proof}
Without loss of generality we may assume that either $\cE$ is a
compact set or $\cE$ contains an infinite interval $[a,\infty)$ for
some $a\in\bbR$. Indeed, if $\cE$ does not contain an infinite
interval then there is an increasing subsequence of gaps
$(a_{j_k},b_{j_k})$ with $b_{j_{k}}<a_{j_{k+1}}$ and
$a_{j_k}\to\infty$ as $k\to\infty$ which splits the set $\cE$ into a
countable disjoint union of compact homogeneous sets
$\cE_0=\cE\cap[E_0,a_{j_{1}}]$, $\cE_k = \cE \cap
[b_{j_{k}},a_{j_{k+1}}]$, $k\in\bbN$. Moreover, it follows from the
proof of \cite[Theorem 2.7]{PV98} that the ratio
$G_{\cE_k}(z,\la_0)/G_{\cE}(z,\la_0)$ of the Green's functions
associated with $\cE_k$ and $\cE$ is bounded in some sufficiently
small neighborhood of $\cE_k$, hence one easily verifies that
$H\in\cR(\cE_k)$ for all $k\geq0$.

Next, fix $x_0\in\bbR$. Then using Lemma \ref{l4.4}, we introduce the
Herglotz function $r(\cdot,x_0)$ by
\begin{align}
r(\ze,x_0) = M_{0,0}(z(\ze),x_0) = d(x_0) +
\int_\bbR\f{d\mu(\la,x_0)}{\la-z}, \quad z\in\bbC_+,
\end{align}
where $M_{0,0}$ is defined in \eqref{4.11a} and
$z(\ze)=(\la_0-\ze)^{-1}$ with $\la_0\in\bbR\bs\si(H)$ introduced in
Hypothesis \ref{h4.3}\,$(iii)$. Abbreviating by $\wti\cE$ and
$\wti\cO$ the preimages of the sets $\cE$ and $\cO$ under the
bijective map $z$,
\begin{align}
\wti\cE=z^{-1}(\cE), \quad \wti\cO=z^{-1}(\cO),
\end{align}
and adding the point zero to $\wti\cE$, and a sufficiently small
neighborhood of zero to $\wti\cO$ in the case of an unbounded set
$\cE$, we note that $\wti\cE$ is compact and homogeneous, $\wti\cO$
is open, $\wti\cE\subset\wti\cO$, and $\ol{\wti\cO}\cap
z^{-1}(\si_\ess(H)\bs\cE) = \emptyset$. Then the functions
$r_j(z,x_0)$, $j=1,2$, defined by
\begin{align}
r_1(z,x_0) &= \int_{\wti\cO} \f{d\mu(\la,x_0)}{\la-z} = \sum_{\la\in
\wti\cO\bs\wti\cE}\f{\mu(\{\la\},x_0)}{\la-z} + \int_{\wti\cE}
\f{d\mu(\la,x_0)}{\la-z}, \quad z\in\bbC\bs\wti\cO,
\lb{4.55} \\
r_2(z,x_0) &= \int_{\bbR\bs\wti\cO} \f{d\mu(\la,x_0)}{\la-z}, \quad
z\in(\bbC\bs\bbR)\cup\wti\cO,
\end{align}
satisfy
\begin{align}
r(z,x_0) = r_1(z,x_0)+r_2(z,x_0), \quad z\in\bbC\bs\bbR.
\end{align}

Since $H\in\cR(\cE)$, one has $\xi(\la,x_0)=1/2$ and hence
$\Re[M_{0,0}(\la,x_0)]=0$ for a.e. $\la\in\cE$. This yields
$\Re[r(\la+i0,x_0)] = 0$ for a.e. $\la\in\wti\cE$, and hence,
\begin{align}
\Re[r_1(\la+i0,x_0)] = -\Re[r_2(\la+i0,x_0)] \,\text{ for a.e. }
\la\in\wti\cE.
\end{align}
Observing that the function $r_2(\cdot,x_0)$ is analytic on
$(\bbC\bs\bbR)\cup\wti\cO$ and $\wti\cE\subset\wti\cO$, one concludes that
$r_2(\cdot,x_0)$ is bounded on $\wti\cE$, and hence,
\begin{align}
\Re[r_1(\cdot+i0,x_0)] = -\Re[r_2(\cdot+i0,x_0)] \in
L^1\big(\wti\cE;dx\big). \lb{4.59}
\end{align}
Moreover, it follows from \cite[Proposition 5.1]{Fi83} that
\eqref{4.21} is equivalent to
\begin{align}
\sum_{\la\in\wti\cO\bs\wti\cE} G_{\wti\cE}(\la,\infty)<\infty.
\lb{4.59a}
\end{align}
In addition, as a consequence of Theorem \ref{tB.7}, the set of mass
points of $d\Om_{0,0}$ is a subset of the set of discrete
eigenvalues of $H$, hence the set of discrete mass points of
$d\mu(\cdot,x_0)$ is a subset of $\wti\cO\bs\wti\cE$. Then it
follows from \eqref{4.55}, \eqref{4.59}, and \eqref{4.59a} that the
function $r_1(\cdot,x_0)$ satisfies the assumptions of Theorem
\ref{tB.9}. Thus, the restriction $d\mu(\cdot,x_0)|_{\wti\cE}$ of
the measure $d\mu(\cdot,x_0)$ to the set $\wti\cE$ is purely
absolutely continuous,
\begin{align}
& d\mu(\cdot,x_0)\big|_{\wti\cE} =
d\mu_\ac(\cdot,x_0)\big|_{\wti\cE}
=\f{1}{\pi}\Im[r_1(\cdot+i0,x_0)]d\la\big|_{\wti\cE},
\end{align}
and hence, it follows from Lemma \ref{l3.4} that
\begin{align}
d\Om_{0,0}(\cdot,x_0)\big|_\cE = d\Om_{0,0,\ac}(\cdot,x_0)\big|_\cE.
\end{align}

Next, performing a similar analysis for the function $M_{1,1}$
defined in \eqref{4.11b}, one obtains
\begin{align}
d\Om_{1,1}(\cdot,x_0)\big|_\cE = d\Om_{1,1,\ac}(\cdot,x_0)\big|_\cE,
\end{align}
and hence,
\begin{align}
d\Om^{\tr}(\cdot,x_0)\big|_\cE =
d\Om^{\tr}_{\ac}(\cdot,x_0)\big|_\cE,
\end{align}
for the restriction of the trace measure $d\Om^{\tr}(\cdot,x_0) =
d\Om_{0,0}(\cdot,x_0) + d\Om_{1,1}(\cdot,x_0)$ associated with $H$.
By \eqref{B.28} and Theorem \ref{tB.7}\,$(i)$ this completes the
proof of \eqref{4.53}.

Finally, equations \eqref{4.11} and \eqref{4.19} imply
\begin{equation}
-1/g(\la+i0,x_0)=\pm 2i \, \Im[m_\pm(\la+i0,x_0)] \, \text{ for a.e.\
$\la\in\cE$.}    \lb{4.30}
\end{equation}
Thus, combining \eqref{4.19}, \eqref{4.30}, and \eqref{B.38} then
yields that the absolutely continuous spectrum of $H$ has uniform
spectral multiplicity two on $\cE$ since
\begin{equation}
\text{for a.e.\ $\la\in\cE$, } \; 0<\pm\Im[m_\pm(\la+i0,x_0)]<\infty.
\end{equation}
\end{proof}
%%%%%%%%%%%%%%%%%%%%%%%%%%%%%%%%%%%%%%%

%%%%%%%%%%%%%%%%%%%%%%%%%%%%%%%%%%%%%%%
%%%%%%%%%%%%% appendices %%%%%%%%%%%%%%
%%%%%%%%%%%%%%%%%%%%%%%%%%%%%%%%%%%%%%%
\appendix
%%%%%%%%%%%%%%%%%%%%%%%%%%%%%%%%%%%%%%%
%%%%%%%%%%%%% appendix B %%%%%%%%%%%%%%
%%%%%%%%%%%%%%%%%%%%%%%%%%%%%%%%%%%%%%%
\section{Herglotz and Weyl--Titchmarsh
Functions \\ for Jacobi and Schr\"odinger Operators in a Nutshell} \lb{B}
\renewcommand{\theequation}{A.\arabic{equation}}
\renewcommand{\thetheorem}{A.\arabic{theorem}}
\setcounter{theorem}{0}
\setcounter{equation}{0}
%%%%%%%%%%%%%%%%%%%%%%%%%%%%%%%%%%%%%%%
%%%%%%%%%%%%%%%%%%%%%%%%%%%%%%%%%%%%%%%

The material in this appendix is known, but since we use it repeatedly at
various places in this paper, we thought it worthwhile to collect it
in an appendix.

%%%%%%%%%%%%%%%%%%%%%%%%%%%%%%%%%%%%%%%
\begin{definition} \lb{dB.1}
Let $\bbC_\pm=\{z\in\bbC \mid \Im(z)\gtrless 0 \}$.
$m:{\mathbb{C_+}}\to {\mathbb{C}}$ is called a {\it Herglotz}
function $($or {\it Nevanlinna} or {\it Pick} function\,$)$ if $m$ is analytic
on
${\mathbb{C}}_+$ and
$m({\mathbb{C}}_+)\subseteq {\mathbb{C}}_+$.
\end{definition}
%%%%%%%%%%%%%%%%%%%%%%%%%%%%%%%%%%%%%%%

\smallskip

\noindent One then extends $m$ to $\bbC_-$ by reflection, that is, one
defines
\begin{equation}
m(z)=\overline{m(\overline z)},
\quad z\in{\mathbb{C}}_-. \lb{B.1}
\end{equation}
Of course, generally, \eqref{B.1} does not represent an analytic
continuation of $m\big|_{\bbC_+}$ into $\bbC_-$.

Fundamental results on Herglotz functions and their representations
on Borel transforms, in part, are due to Fatou, Herglotz, Luzin,
Nevanlinna, Plessner, Privalov, de la Vall{\'e}e Poussin, Riesz, and
others.  Here we just summarize a few of these results:

%%%%%%%%%%%%%%%%%%%%%%%%%%%%%%%%%%%%%%
\begin{theorem} $($\cite[Sect.\ 69]{AG81}, \cite{AD56},
\cite[Chs.\ II, IV]{Do74}, \cite{KK74}, \cite[Ch.\ 6]{Ko98}, \cite[Chs.\ II, IV]{Pr56},
 \cite[Ch.\ 5]{RR94}$)$. \lb{tB.2} Let $m$ be a Herglotz function. Then, \\
$(i)$ There exists a nonnegative measure $d\om$ on ${\mathbb{R}}$
satisfying
\begin{equation} \lb{B.2}
\int_{{\mathbb{R}}} \frac{d\om (\la )}{1+\la^2}<\infty
\end{equation}
such that the Nevanlinna, respectively, Riesz-Herglotz
representation
\begin{align}
\begin{split}
&m(z)=c+dz+\int_{{\mathbb{R}}} d\om (\la) \bigg(\frac{1}{\la
-z}-\frac{\la}
{1+\la^2}\bigg), \quad z\in\bbC_+, \lb{B.3} \\[2mm]
& \, c=\Re[m(i)],\quad d=\lim_{\eta \uparrow
\infty}m(i\eta )/(i\eta ) \geq 0
\end{split}
\end{align}
holds. Conversely, any function $m$ of the type \eqref{B.3} is a
Herglotz function. \\
$(ii)$ The absolutely continuous $({\it ac})$ part $d\om_{ac}$ of
$d\om$ with respect to Lebesgue measure $d\la$ on ${\mathbb{R}}$ is
given by
\begin{equation}\lb{B.5}
d\om_{ac}(\la)=\pi^{-1}\Im[m(\la+i0)]\,d\la.
\end{equation}
\end{theorem}
%%%%%%%%%%%%%%%%%%%%%%%%%%%%%%%%%%%%%

Next, we denote by
\begin{equation}
d\mu =d\mu_{\ac}+d\mu_{\sc} +d\mu_{\pp} \lb{B.8}
\end{equation}
the decomposition of a measure $d\mu$ into its absolutely continuous
$(\ac)$, singularly continuous $(\sc)$, and pure point $(\pp)$ parts
with respect to Lebesgue measure on ${\mathbb{R}}$.

%%%%%%%%%%%%%%%%%%%%%%%%%%%%%%%%%%%%
\begin{theorem} $($\cite{Gi84}, \cite{GP87}$)$.  \lb{tB.5}
Let $m$ be a Herglotz function with representation \eqref{B.3} and
denote by $\La$ the set
\begin{equation}
\La=\{\la\in\bbR\,|\, \Im[m(\la+i0)] \, \text{exists $($finitely or
infinitely$)$}\}. \lb{B.9}
\end{equation}
Then, $S$, $S_{\ac}$, $S_{\rm s}$, $S_{\sc}$, $S_{\pp}$ are
essential supports of $d\om$, $d\om_{\ac}$, $d\om_{\rm s}$,
$d\om_{\sc}$, $d\om_{\pp}$, respectively, where
\begin{align}
S&=\{\la\in\La\,|\, 0<\Im[m(\la+i0)]\leq\infty\}, \lb{B.10}
\\
S_{\ac}&=\{\la\in\La\,|\, 0<\Im[m(\la+i0)]<\infty\},
\lb{B.11} \\
S_{\rm s}&=\{\la\in\La\,|\, \Im[m(\la+i0)]=\infty\},
\lb{B.12} \\
S_{\sc}&=\Big\{\la\in\La\,|\, \Im[m(\la+i0)]=\infty, \,
\lim_{\eps\downarrow 0}
(-i\eps)m(\la+i\eps)=0\Big\}, \lb{B.13} \\
S_{\pp}&=\Big\{\la\in\La\,|\, \Im[m(\la+i0)]=\infty, \,
\lim_{\eps\downarrow 0} (-i\eps)m(\la+i\eps)=\om(\{\la\})>0\Big\}.
\lb{B.14}
\end{align}
\end{theorem}
%%%%%%%%%%%%%%%%%%%%%%%%%%%%%%%%%%%%%

Next, consider Herglotz functions $\pm m_\pm$ of the type
\eqref{B.3},
\begin{align}
\begin{split}
& \pm m_\pm(z)=c_\pm+d_\pm z+\int_{{\mathbb{R}}} d\om_\pm (\la)
\bigg(\frac{1}{\la -z}-\frac{\la}
{1+\la^2}\bigg), \quad z\in\bbC_+, \lb{B.21} \\
& \,\, c_\pm\in\bbR, \quad d_\pm\geq 0,
\end{split}
\end{align}
and introduce the $2\times 2$ matrix-valued Herglotz function $M$
\begin{align}
&M(z)=\big(M_{j,k}(z)\big)_{j,k=0,1}, \quad z\in\bbC_+, \lb{B.22} \\
&M(z)=\f{1}{m_-(z)-m_+(z)}\begin{pmatrix} 1 &
\f{1}{2}[m_-(z)+m_+(z)]
\\
\f{1}{2}[m_-(z)+m_+(z)] & m_-(z)m_+(z) \end{pmatrix}
\lb{B.23} \\
& \hspace*{.85cm} =C+Dz + \int_{{\mathbb{R}}} d\Om (\la)
\bigg(\frac{1}{\la -z}-\frac{\la}
{1+\la^2}\bigg), \quad z\in\bbC_+,  \lb{B.24} \\
& \, C=C^*, \quad D \geq 0 \no
\end{align}
with $C=(C_{j,k})_{j,k=0,1}$ and $D=(D_{j,k})_{j,k=0,1}$ $2\times 2$
matrices and $d\Om=(d\Om_{j,k})_{j,k=0,1}$ a $2\times 2$
matrix-valued nonnegative measure satisfying
\begin{equation}
\int_\bbR \f{d|\Om_{j,k}(\la)|}{1+\la^2} < \infty, \quad j,k=0,1.
\lb{B.25}
\end{equation}
Moreover, we introduce the trace Herglotz function $M^{\tr}$
\begin{align}
&M^{\tr}(z)=M_{0,0}(z)+M_{1,1}(z)=\f{1+m_-(z)m_+(z)}{m_-(z)-m_+(z)}
\lb{B.26} \\
& \hspace*{1.05cm} =c + dz + \int_{\bbR} d\Om^{\tr}(\la)
\bigg(\frac{1}{\la -z}-\frac{\la} {1+\la^2}\bigg), \quad z\in\bbC_+,
\lb{B.27} \\ & \, c\in\bbR, \;  d\geq 0, \quad d\Om^{\tr}=
d\Om_{0,0} + d\Om_{1,1}. \no
\end{align}
Then,
\begin{equation}
d\Om \ll d\Om^{\tr} \ll d\Om   \lb{B.28}
\end{equation}
(where $d\mu \ll d\nu$ denotes that $d\mu$ is absolutely continuous
with respect to $d\nu$).

The next result holds both for the Jacobi and Schr\"odinger cases.
In the Jacobi case we identify
\begin{equation}
m_\pm(z) \, \text{ and } \, M_\pm(z,n_0), \quad z\in\bbC_+,
\end{equation}
where $M_\pm(z,n_0)$ denote the half-lattice Weyl--Titchmarsh
$m$-functions defined in \eqref{2.10}--\eqref{2.12} and in the
Schr\"odinger case
\begin{equation}
m_\pm(z) \, \text{ and } \, m_\pm(z,x_0), \quad z\in\bbC_+,
\end{equation}
where $m_\pm(z,x_0)$ are the half-line Weyl--Titchmarsh
$m$-functions defined in \eqref{4.10}. One then has the following
basic result.

%%%%%%%%%%%%%%%%%%%%%%%%%%%%%%%%%%%%
\begin{theorem} $($\cite{Gi98}, \cite{Ka62}, \cite{Ka63}, \cite{Si05b},
\cite{Te00}$)$.  \lb{tB.7} ${}$ \\
$(i)$ The operator $H$ $($in the Jacobi case $H$ is defined in
\eqref{2.4} and in the Schr\"odinger case in \eqref{4.4}$)$ is
unitarily equivalent to the operator of multiplication by $I_2\id$
$($where $I_2$ is the $2\times 2$ identity matrix and
$\id(\la)=\la$, $\la\in\bbR$$)$ on $L^2(\bbR; d\Om(\cdot))$, and
hence,
\begin{align}
\si(H) = \supp\,(d\Om) = \supp\,(d\Om^\tr),
\end{align}
where $d\Om$ and $d\Om^\tr$ are introduced in \eqref{B.24} and
\eqref{B.27}, respectively. \\
$(ii)$ The spectral multiplicity of $H$ is two if and only if
\begin{equation}
|\cM_2|>0,  \lb{B.37}
\end{equation}
where
\begin{equation}
\cM_2=\{\la\in\La_+\,|\,
m_+(\la+i0)\in\bbC\bs\bbR\}\cap\{\la\in\La_-\,|\,
m_-(\la+i0)\in\bbC\bs\bbR\}.  \lb{B.38}
\end{equation}
If $|\cM_2|=0$, the spectrum of $H$ is simple. Moreover, $\cM_2$ is
a maximal set on which $H$ has uniform multiplicity two. \\
$(iii)$ A maximal set $\cM_1$ on which $H$ has uniform multiplicity
one is given by
\begin{align}
\cM_1&=\{\la\in\La_+\cap\La_-\,|\, m_+(\la+i0)=
m_-(\la+i0)\in\bbR\}  \no \\
&\quad \cup \{\la\in\La_+\cap\La_-\,|\,
|m_+(\la+i0)|= |m_-(\la+i0)|=\infty\}  \no \\
&\quad \cup \{\la\in\La_+\cap\La_-\,|\, m_+(\la+i0)\in\bbR,
m_-(\la+i0)\in\bbC\bs\bbR\}  \no \\
&\quad \cup \{\la\in\La_+\cap\La_-\,|\, m_-(\la+i0)\in\bbR,
m_+(\la+i0)\in\bbC\bs\bbR\}. \lb{B.39}
\end{align}
In particular, $\sigma_{\rm s}(H)=\sigma_{\sc}(H)\cup
\sigma_{\pp}(H)$ is always simple.
\end{theorem}
%%%%%%%%%%%%%%%%%%%%%%%%%%%%%%%%%%%%%

Finally, we give a proof of a basic result due to Peherstorfer and
Yuditskii \cite[Lemma 2.4]{PY03}.

%%%%%%%%%%%%%%%%%%%%%%%%%%%%%%%%%%%%%%%
\begin{theorem} $($\cite[Lemma 2.4]{PY03}$)$.   \lb{tB.9}
Let $\cE\subset\bbR$ be a compact homogenous set and $m(z)$ a
Herglotz function with the representation
\begin{align}
\begin{split}
& m(z)=a + \sum_{j\in J}\frac{\mu(\{\la_j\})}{\la_j-z} +
\int_{\cE} \f{d\mu(\la)}{\la-z}, \quad z\in\bbC_+, \\
& a\in\bbR, \;\; J\subseteq\bbN, \;\; d\mu \text{ a finite measure},
\;\; \supp \, (d\mu) \subseteq \cE\cup\{\la_j\}_{j\in J}.
\end{split}
\end{align}
Denote by $m(\la+i0)=\lim_{\eps\downarrow 0}m(\la+i\eps)$ the a.e.\
normal boundary values of $m$ and assume that
\begin{align}
&\Re[m(\cdot+i0)] \in L^1\big(\cE;d\la\big), \lb{B.48} \\
&\sum_{j\in J} G_{\cE}(\la_j,\infty)<\infty, \lb{B.49}
\end{align}
where $G_{\cE}(\cdot,\infty)$ is the potential  theoretic Green's
function for the domain $(\bbC\cup\{\infty\})\bs\cE$ with
logarithmic singularity at infinity $($cf., e.g., \cite[Sect.\ 5.2]{Ra95}$)$,
\begin{align}
G_\cE(z,\infty) \underset{|z|\to\infty}{=} \log |z| - \log (\ca(\cE))
+ \oh(1). \lb{B.50}
\end{align}
Then, $d\mu$ is purely absolutely continuous and hence
\begin{equation}
d\mu\big|_\cE=d\mu_{\ac}\big|_\cE=\f{1}{\pi}
\Im[m(\cdot+i0)]\,d\la\big|_\cE.
\end{equation}
\end{theorem}
%%%%%%%%%%%%%%%%%%%%%%%%%%%%%%%%%%%%%%%
\begin{proof}
First, we briefly introduce some important notation (cf.,
\cite{CSZ1}, \cite{PY03}, \cite[Ch.\ 8]{Si08a}, \cite{SY97} for
a comprehensive discussion). Let $\Ga$ be the Fuchsian group of
linear-fractional transformations of the unit disc $\bbD$,
uniformizing the domain $\Om=(\bbC\cup\{\infty\})\bs\cE$, and denote
by $\Ga^*$ the group of unimodular characters associated with the
Fuchsian group $\Ga$ (i.e., each $\al\in\Ga^*$ is a homomorphism
from $\Ga$ to $\dD$). By
\begin{align}
\bbF = \{z\in\bbD \st |\ga'(z)|<1 \text{ for all }
\ga\in\Ga\bs\{\id\}\}
\end{align}
we denote the Ford orthocircular fundamental domain of $\Ga$ (see
\cite{Fo51}, \cite{SY97}), and as usual, keeping the same notation,
we add to $\bbF$ half of its boundary circles (say, the boundary
circles lying in $\bbC_+\cap\bbD$). By $\cm(z)$ we denote the
uniformizing map (also known as the universal covering map for
$\Om$), that is, $\cm(z)$ is the unique map satisfying the following
properties:
\begin{enumerate}[$(i)$]
\item
$\cm(z)$ maps $\bbD$ meromorphically onto $\Om$.
\item
$\cm(z)$ is $\Ga$-automorphic, that is, $\cm\circ\ga = \cm$,
$\ga\in\Ga$.
\item
$\cm(z)$ is locally a bijection ($\cm$ maps $\bbF$ bijectively to
$\Om$).
\item
$\cm(0)=\infty$ and $\lim_{z\to0}z\cm(z)>0$.
\end{enumerate}
We also introduce the notion of Blaschke products associated with
the group $\Ga$,
\begin{align}
B(z,w) = \prod_{\ga\in\Ga}
\f{\ga(w)-z}{1-\ol{\ga(w)}z}\frac{|\ga(w)|}{\ga(w)}, \quad z\in\bbD,
\end{align}
where we set $\frac{|\ga(w)|}{\ga(w)}\equiv-1$ if $\ga(w)=0$. Then
we define
\begin{align}
B(z) = B(z,0), \quad B_\infty(z)=\prod_{j\in J}B(z,z_j), \quad
z\in\bbD,
\end{align}
(condition \eqref{B.49} guarantees convergence of the above
product), where $\{z_j\}_{j\in J}\subset\bbF$ are the points
satisfying $\cm(z_j)=\la_j$, $j\in J$. It follows that the functions
$B$ and $B_\infty$ are character-automorphic (i.e., there are
characters $\al,\be\in\Ga^*$ such that $B\circ\ga=\al(\ga)B$ and
$B_\infty\circ\ga=\be(\ga)B_\infty$ for all $\ga\in\Ga$) and the
following formulas hold
\begin{align}
&G_\cE(\cm(z),\infty)=-\log\big(|B(z)|\big), \\
&\lim_{z\to0}B(z)\cm(z)=\ca(\cE), \quad
\lim_{z\to0}\f{B(z)^2\cm'(z)}{B'(z)}=-\ca(\cE), \lb{B.55}
\end{align}
where $\ca(\cE)$ denotes logarithmic capacity of the set $\cE$ (cf.\
\eqref{B.50}). Moreover, we define
\begin{align}
\phi(z)=\frac{zB'(z)}{B(z)}, \quad z\in\bbD,
\end{align}
then one verifies
\begin{align}
\phi(\ze) = \sum_{\ga\in\Ga}|\ga'(\ze)| \,\text{ for a.e.
}\,\ze\in\dD,
\end{align}
and
\begin{align}
\phi(\ga(\ze))|\ga'(\ze)|=\phi(\ze) \,\text{ for all }\, \ga\in\Ga,
\; \ze\in\dD. \lb{B.56}
\end{align}

After these preliminaries we commence with the proof of Theorem
\ref{tB.9}. One observes that it suffices to prove
\begin{align}
\mu(\bbR) = \sum_{j\in J}\mu(\{\la_j\}) + \f{1}{\pi}\int_\cE
\Im[m(\la+i0)]\,d\la. \lb{B.57}
\end{align}
Let $r(z)$ denote the Herglotz function
\begin{align}
r(z)=m(z)-a=\sum_{j\in J}\frac{\mu(\{\la_j\})}{\la_j-z} + \int_{\cE}
\f{d\mu(\la)}{\la-z}, \quad z\in\bbC_+. \lb{B.58}
\end{align}
Then one verifies the following asymptotic formula,
\begin{align}
r(z) \underset{z\to\infty}{=} -\frac{\mu(\bbR)}{z} + \Oh(z^{-2}).
\lb{B.59}
\end{align}
Using the symmetry $r(\ol{z})=\ol{r(z)}$ and \eqref{B.58} one also
derives
\begin{align}
r(\la+i0)-r(\la-i0) = 2i\Im[r(\la+i0)] = 2i\Im[m(\la+i0)]\,\text{
for a.e. } \la\in\bbR.
\end{align}
Moreover, it follows from condition \eqref{B.48} that $r(\cdot\pm
i0)\in L^1\big(\cE;d\la\big)$.

Next, one computes,
\begin{align}
&\f{1}{\pi}\int_\cE\Im[m(\la+i0)]\,d\la = \f{1}{2\pi
i}\int_\cE[r(\la+i0)-r(\la-i0)]\,d\la = \f{-1}{2\pi i}\oint_\dOm
r(\la+i0)\,d\la \no
\\
&\quad = \f{1}{2\pi i}\int_{\dF\cap\dD} r(\cm(\ze))\cm'(\ze)\,d\ze =
\f{1}{2\pi i}\oint_\dD r(\cm(\ze))\f{\cm'(\ze)}{\phi(\ze)}\,d\ze.
\lb{B.61}
\end{align}
To evaluate the last integral one utilizes the Direct Cauchy Theorem
(cf.\ \cite[Lemma 1.1]{PY03}, \cite[Theorem H]{SY97}),
\begin{align}
&\f{1}{2\pi i}\oint_\dD r(\cm(\ze))\f{\cm'(\ze)}{\phi(\ze)}\,d\ze =
\f{1}{2\pi i}\oint_\dD \f{\f{r(\cm(\ze))B_\infty(\ze)}{B(\ze)}
\f{B(\ze)^2\cm'(\ze)}{\phi(\ze)}}{B(\ze)B_\infty(\ze)}\,d\ze \no
\\
&\quad = \f{r(\cm(0))B(0)^2\cm'(0)}{B(0)B'(0)} + \sum_{j\in J}
\f{r(\cm(z_j))B_\infty(z_j)\cm'(z_j)}{B_\infty'(z_j)}\no
\\
&\quad = \mu(\bbR) - \sum_{j\in J} \mu(\{\la_j\}). \lb{B.65}
\end{align}
Thus, combining \eqref{B.61} and \eqref{B.65} yields \eqref{B.57}.
\end{proof}
%%%%%%%%%%%%%%%%%%%%%%%%%%%%%%%%%%%%%%%

%%%%%%%%%%%%%%%%%%%%%%%%%%%%%%%%%%%%%%%
%%%%%%%%%%%%% appendix B %%%%%%%%%%%%%%
%%%%%%%%%%%%%%%%%%%%%%%%%%%%%%%%%%%%%%%
\section{Caratheodory and Weyl--Titchmarsh Functions
\\ for CMV Operators in a Nutshell} \lb{C}
\renewcommand{\theequation}{B.\arabic{equation}}
\renewcommand{\thetheorem}{B.\arabic{theorem}}
\setcounter{theorem}{0} \setcounter{equation}{0}
%%%%%%%%%%%%%%%%%%%%%%%%%%%%%%%%%%%%%%%
%%%%%%%%%%%%%%%%%%%%%%%%%%%%%%%%%%%%%%%

In this appendix we provide some basic facts on Caratheodory functions
and prove the analog of Theorem \ref{tB.7} for CMV operators.

%%%%%%%%%%%%%%%%%%%%%%%%%%%%%%%%%%%%%%%
\begin{definition} \lb{dC.1}
Let $\bbD$ and $\dD$ denote the open unit disk and the
counterclockwise oriented unit circle in the complex plane $\bbC$,
\begin{equation}
\bbD = \{ z\in\bbC \st \abs{z} < 1 \}, \quad \dD = \{ \ze\in\bbC \st
\abs{\ze} = 1 \},
\end{equation}
and $\Cl$ and $\Cr$ the open left and right complex half-planes,
respectively,
\begin{equation}
\Cl = \{z\in\bbC \st \Re(z) < 0\}, \quad \Cr = \{z\in\bbC \st \Re(z)
> 0\}.
\end{equation}
A function $f:\bbD\to\bbC$ is called Caratheodory if $f$ is analytic
on $\bbD$ and $f(\bbD)\subset\Cr$. One then extends $f$ to
$\bbC\bs\ol{\bbD}$ by reflection, that is, one defines
\begin{align}
f(z) = -\ol{f(1/\ol{z})}, \quad z\in\bbC\bs\ol{\bbD}. \lb{C.1}
\end{align}
Of course, generally, \eqref{C.1} does not represent the analytic
continuation of $f|_\bbD$ into $\bbC\bs\ol{\bbD}$.
\end{definition}
%%%%%%%%%%%%%%%%%%%%%%%%%%%%%%%%%%%%%%%

The fundamental result on Caratheodory functions reads as follows.

%%%%%%%%%%%%%%%%%%%%%%%%%%%%%%%%%%%%%%%
\begin{theorem} $($\cite[Sect.\ 3.1]{Ak65}, \cite[Sect.\ 69]{AG81},
\cite[Sect.\ 1.3]{Si05}$)$.\  \label{tC.2}
Let $f$ be a Caratheodory function. Then, \\
$(i)$ $f(z)$ has finite normal limits $f(\ze)=\lim_{r\uparrow1}
f(r\ze)$ for a.e.~$\ze\in\dD$. \\
$(ii)$ Suppose $f(r\ze)$ has a zero normal limit on a subset of
$\dD$ having positive Lebesgue measure. Then $f\equiv 0$. \\
$(iii)$ There exists a nonnegative finite measure $d\om$ on $\dD$
such that the Herglotz representation
\begin{align}
\begin{split}
& f(z)=ic+ \oint_{\dD} d\om(\zeta) \, \f{\zeta+z}{\zeta-z}, \quad
z\in\bbD,
\\
& c=\Im(f(0)), \quad \oint_{\dD} d\om(\zeta) = \Re(f(0)) < \infty,
\end{split} \lb{C.3}
\end{align}
holds. Conversely, any function $f$ of the type \eqref{C.3} is a
Caratheodory function. \\
$(iv)$ The absolutely continuous $({\it ac})$ part $d\om_{ac}$ of
$d\om$ with respect to the normalized Lebesgue measure $d\om_0$ on
$\dD$ is given by
\begin{equation}\lb{C.5}
d\om_{ac}(\ze)=\pi^{-1}\Re[f(\ze)]\,d\om_0(\ze).
\end{equation}
\end{theorem}
%%%%%%%%%%%%%%%%%%%%%%%%%%%%%%%%%%%%%%%

Next, we denote by
\begin{equation}
d\om =d\om_{\ac}+d\om_{\sc} +d\om_{\pp} \lb{C.8}
\end{equation}
the decomposition of $d\om$ into its absolutely continuous $({\it
ac})$, singularly continuous $({\it sc})$, and pure point $({\it
pp})$ parts with respect to Lebesgue measure on $\dD$.

%%%%%%%%%%%%%%%%%%%%%%%%%%%%%%%%%%%%%%%
\begin{theorem} $($\cite[Sects.\ 1.3, 1.4]{Si05}$)$.  \lb{tC.5}
Let $f$ be a Caratheodory function with representation \eqref{C.3}
and denote by $\La$ the set
\begin{equation}
\La=\{\ze\in\dD \st \Re[f(\ze)] \, \text{exists $($finitely or
infinitely$)$}\}. \lb{C.9}
\end{equation}
Then, $S$, $S_{\ac}$, $S_{\rm s}$, $S_{\sc}$, $S_{\pp}$ are essential
supports of $d\om$, $d\om_{\ac}$, $d\om_{\rm s}$, $d\om_{\sc}$,
$d\om_{\pp}$, respectively, where
\begin{align}
S&=\{\ze\in\La \st 0<\Re[f(\ze)]\leq\infty\}, \lb{C.10}
\\
S_{\ac}&=\{\ze\in\La \st 0<\Re[f(\ze)]<\infty\},
\lb{C.11} \\
S_{\rm s}&=\{\ze\in\La \st \Re[f(\ze)]=\infty\},
\lb{C.12} \\
S_{\sc}&=\Big\{\ze\in\La \st \Re[f(\ze)]=\infty, \,
\lim_{r\uparrow1}(1-r)f(r\ze)=0\Big\}, \lb{C.13} \\
S_{\pp}&=\Big\{\ze\in\La \st \Re[f(\ze)]=\infty, \,
\lim_{r\uparrow1}\left(\f{1-r}{2}\right)f(r\ze)=\om(\{\ze\})>0\Big\}.
\lb{C.14}
\end{align}
\end{theorem}
%%%%%%%%%%%%%%%%%%%%%%%%%%%%%%%%%%%%%%%

Next, consider Caratheodory functions $\pm m_\pm$ of the type
\eqref{C.3},
\begin{align}
\begin{split}
& \pm m_\pm(z)=ic_\pm + \oint_{\dD} d\om_\pm(\zeta) \,
\f{\zeta+z}{\zeta-z},
\quad z\in\bbD, \lb{C.21} \\
& \,\, c_\pm\in\bbR,
\end{split}
\end{align}
and introduce the $2\times 2$ matrix-valued Caratheodory function
$\wti M$ by
\begin{align}
&\wti M(z)=\big(\wti M_{j,k}(z)\big)_{j,k=0,1}, \quad z\in\bbD, \lb{C.22} \\
&\wti M(z)=\f{1}{m_+(z)-m_-(z)}\begin{pmatrix} 1 &
\f{1}{2}[m_+(z)+m_-(z)]
\\
-\f{1}{2}[m_+(z)+m_-(z)] & -m_+(z)m_-(z) \end{pmatrix},
\lb{C.23} \\
& \hspace*{.85cm} = i\wti C + \oint_{{\dD}} d\wti\Om(\ze)
\frac{\ze+z}{\ze-z}, \quad z\in\bbD,  \lb{C.24} \\
& \, \wti C = \wti C^* = \Im[\wti M(0)],\no
\end{align}
where $d\wti\Om=(d\wti\Om_{j,k})_{j,k=0,1}$ is a $2\times 2$
matrix-valued nonnegative measure satisfying
\begin{equation}
\oint_\dD d\,|\wti\Om_{j,k}(\ze)| < \infty, \quad j,k=0,1. \lb{C.25}
\end{equation}
Moreover, we introduce the trace Caratheodory function $\wti
M^{\tr}$
\begin{align}
&\wti M^{\tr}(z)=\wti M_{0,0}(z) + \wti M_{1,1}(z) =
\f{1-m_+(z)m_-(z)}{m_+(z)-m_-(z)}
\lb{C.26} \\
& \hspace*{1.05cm} =i\wti c + \oint_{\dD} d\wti\Om^{\tr}(\ze)
\frac{\ze+z}{\ze-z}, \quad z\in\bbD, \lb{C.27} \\
& \, \wti c\in\bbR, \quad d\wti\Om^{\tr}= d\wti\Om_{0,0} +
d\wti\Om_{1,1}. \no
\end{align}
Then,
\begin{equation}
d\wti\Om \ll d\wti\Om^{\tr} \ll d\wti\Om   \lb{C.28}
\end{equation}
(where $d\mu \ll d\nu$ denotes that $d\mu$ is absolutely continuous
with respect to $d\nu$). This implies that there is a self-adjoint
integrable $2\times 2$ matrix $\wti R(\ze)$ such
\begin{align}
d\wti\Om(\ze) = \wti R(\ze) d\wti\Om^\tr(\ze) \lb{C.29}
\end{align}
by the Radon--Nikodym theorem. Moreover, the matrix $\wti R(\ze)$ is
nonnegative and given by
\begin{align}
\wti R(\ze) = \lim_{r\uparrow1}\frac{1} {\Re\big[\wti M_{0,0}(r\ze)+\wti
M_{1,1}(r\ze)\big]}
\begin{pmatrix}
\Re\big[\wti M_{0,0}(r\ze)\big] & i\Im\big[\wti M_{0,1}(r\ze)\big] \\
i\Im\big[\wti M_{1,0}(r\ze)\big] & \Re\big[\wti M_{1,1}(r\ze)\big]
\end{pmatrix} \no
\\
\text{ for a.e.\ $\ze\in\dD$}. \lb{C.30}
\end{align}

Next, we identify
\begin{equation}
m_\pm(z) \, \text{ and } \, M_\pm(z,n_0), \quad n_0\in\bbZ,\;
z\in\bbD, \lb{C.31}
\end{equation}
where $M_\pm(z,n_0)$ denote the half-lattice Weyl--Titchmarsh
$m$-functions defined in \eqref{3.10}--\eqref{3.11} and introduce
another $2\times 2$ matrix-valued Caratheodory function
\begin{align}
& M(z,n_0) =
\begin{pmatrix}
M_{0,0}(z,n_0) & M_{0,1}(z,n_0) \\
M_{1,0}(z,n_0) & M_{1,1}(z,n_0)
\end{pmatrix} \no
\\
&\quad = \begin{pmatrix} (\de_{n_0-1},(U+zI)(U-zI)^{-1}\de_{n_0-1})
&(\de_{n_0-1},(U+zI)(U-zI)^{-1}\de_{n_0})
\\
(\de_{n_0},(U+zI)(U-zI)^{-1}\de_{n_0-1}) &
(\de_{n_0},(U+zI)(U-zI)^{-1}\de_{n_0})
\end{pmatrix} \no
\\
&\quad = \oint_\dD d\Omega(\ze,n_0)\, \frac{\ze+z}{\ze-z}, \quad
z\in\bbD, \lb{C.32}
\end{align}
where $U$ denotes a CMV operator of the form \eqref{3.4} and
$d\Om=(d\Om_{j,k})_{j,k=0,1}$ a $2\times 2$ matrix-valued
nonnegative measure satisfying
\begin{equation}
\oint_\dD d\,|\Om_{j,k}(\ze)| < \infty, \quad j,k=0,1.
\end{equation}
Then the two Caratheodory matrices $M(z,n_0)$ and $\wti M(z)$ are
related by (cf.\ \cite[Equation (3.62)]{GZ06})
\begin{align}
&\wti M(z) +
\begin{pmatrix}
\f{i}{2}\Im(\alpha_{n_0}) & \f{1}{2}\Re(\alpha_{n_0}) \\
-\f{1}{2}\Re(\alpha_{n_0}) & -\f{i}{2}\Im(\alpha_{n_0})
\end{pmatrix} \no
\\
&\quad =\begin{cases}
\f{1}{4}\begin{pmatrix} \rho_{n_0} & \rho_{n_0} \\
-b_{n_0} & a_{n_0} \end{pmatrix}^* M(z,n_0) \begin{pmatrix}
\rho_{n_0} & \rho_{n_0}
\\[1mm]
-b_{n_0} & a_{n_0} \end{pmatrix}, & \text{$n_0$ odd}, \\
\f{1}{4}\begin{pmatrix} -\rho_{n_0} & \rho_{n_0} \\
\ol{b_{n_0}} & \ol{a_{n_0}} \end{pmatrix}^* M(z,n_0)
\begin{pmatrix} -\rho_{n_0} &
\rho_{n_0} \\
\ol{b_{n_0}} & \ol{a_{n_0}} \end{pmatrix}, & \text{$n_0$ even},
\end{cases}, \quad z\in\bbD.
\end{align}
Moreover, it follows from \cite[Lemma 3.2]{GZ06} that
\begin{align}
M_{1,1}(z,n_0)=\frac{1-M_+(z,n_0)M_-(z,n_0)}{M_+(z,n_0)-M_-(z,n_0)},
\quad z\in\bbD, \lb{C.34}
\end{align}
hence by \eqref{C.26} and \eqref{C.31}
\begin{align}
M_{1,1}(z,n_0)=\wti M^\tr(z), \quad z\in\bbD. \lb{C.35}
\end{align}

One then has the following basic result (see also \cite{Si05b}).

%%%%%%%%%%%%%%%%%%%%%%%%%%%%%%%%%%%%%%%
\begin{theorem}\lb{tC.7} ${}$ \\
$(i)$ The CMV operator $U$ on $\ltz$ defined in \eqref{3.4} is unitarily
equivalent to the operator of multiplication by $I_2\id$ $($where
$I_2$ is the $2\times 2$ identity matrix and $\id(\ze)=\ze$,
$\ze\in\dD$$)$ on $L^2(\dD; d\wti\Om(\cdot))$, and hence,
\begin{align}
\si(U) = \supp\,(d\wti\Om) = \supp\,(d\wti\Om^\tr), \lb{C.36}
\end{align}
where $d\wti\Om$ and $d\wti\Om^\tr$ are introduced in \eqref{C.24}
and
\eqref{C.27}, respectively. \\
$(i')$ The operator $U$ is also unitarily equivalent to the operator
of multiplication by $I_2\id$ on $L^2(\dD; d\Om(\cdot))$, and hence
by $(i)$, \eqref{C.28}, and \eqref{C.35},
\begin{align}
d\Om \ll d\Om_{1,1} \ll d\Om  \,\text{ and }\, \si(U) =
\supp\,(d\Om) = \supp\,(d\Om_{1,1}), \lb{C.36a}
\end{align}
where $d\Om$ and $d\Om_{1,1}$ are introduced in \eqref{C.32}. \\
$(ii)$ The spectral multiplicity of $U$ is two if and only if
\begin{equation}
|\cM_2|>0,  \lb{C.37}
\end{equation}
where
\begin{equation}
\cM_2=\{\ze\in\La_+\,|\, m_+(\ze)\in\bbC\bs
i\bbR\}\cap\{\ze\in\La_-\,|\, m_-(\ze)\in\bbC\bs i\bbR\}.  \lb{C.38}
\end{equation}
If $|\cM_2|=0$, the spectrum of $U$ is simple. Moreover, $\cM_2$ is
a maximal set on which $U$ has uniform multiplicity two. \\
$(iii)$ A maximal set $\cM_1$ on which $U$ has uniform multiplicity
one is given by
\begin{align}
\cM_1&=\{\ze\in\La_+\cap\La_-\,|\, m_+(\ze)=
m_-(\ze)\in i\bbR\}  \no \\
&\quad \cup \{\ze\in\La_+\cap\La_-\,|\,
|m_+(\ze)|= |m_-(\ze)|=\infty\}  \no \\
&\quad \cup \{\ze\in\La_+\cap\La_-\,|\, m_+(\ze)\in i\bbR,
m_-(\ze)\in\bbC\bs i\bbR\}  \no \\
&\quad \cup \{\ze\in\La_+\cap\La_-\,|\, m_-(\ze)\in i\bbR,
m_+(\ze)\in\bbC\bs i\bbR\}. \lb{C.39}
\end{align}
In particular, $\sigma_{\rm s}(U)=\sigma_{\sc}(U)\cup
\sigma_{\pp}(U)$ is always simple.
\end{theorem}
%%%%%%%%%%%%%%%%%%%%%%%%%%%%%%%%%%%%%%%
\begin{proof}
We refer to Lemma 3.6 and Corollary 3.5 in \cite{GZ06} for a proof
of $(i)$ and $(i')$, respectively. To prove $(ii)$ and $(iii)$
observe that by $(i)$ and \eqref{C.29}
\begin{align}
\cN_k = \{\ze\in\si(U) \st \rank[\wti R(\ze)]=k\}, \quad k=1,2,
\end{align}
denote the maximal sets where the spectrum of $U$ has multiplicity
one and two, respectively. Using \eqref{C.23} and \eqref{C.30} one
verifies that $\cN_k = \cM_k$, $k=1,2$.
\end{proof}
%%%%%%%%%%%%%%%%%%%%%%%%%%%%%%%%%%%%%%%

\medskip
%%%%%%%%%%%%%%%%%%%%%
\noindent {\bf Acknowledgments.} We are indebted to Barry Simon and
Alexander Volberg for helpful discussions on this topic.
%%%%%%%%%%%%%%%%%%

%%%%%%%%%%%%%%%%%%%%%%%%%%%%%%%%%%%%%%
%%%%%%%%%%%%%%%%%%%%%%%%%%%%%%%%%%%%%%

\end{document}